\documentclass[11pt]{article}
\usepackage{latexsym,amssymb,amsmath,amscd,amsthm,amsxtra}
\usepackage{mathrsfs}
\usepackage{epsfig}
\usepackage{pgf,tikz}
\usepackage{graphicx}
\usepackage{graphics}
\usepackage{setspace}
\newcommand{\counte}{section}

\newtheorem{lemma}{\bf Lemma}[\counte]

\author{Zhao Xu-an, zhaoxa@bnu.edu.cn\\ Department of Mathematics, Beijing Normal University\\Key Laboratory
of Mathematics and Complex Systems\\ Ministry of Education,
China, Beijing 100875\\Jin Chunhua, jinch@amss.ac.cn\\Academy of Mathematics and Systems Science, Chinese Academy of Sciences
\\China, Beijing 100190}

\title{The ranks of homotopy groups of Kac-Moody groups\thanks{The authors are supported by NSFC 11171025}}
\date{}

\begin{document}

\maketitle
\begin{abstract}
Let $A$ be a Cartan matrix and $G(A)$ be the Kac-Moody group associated to Cartan matrix $A$. In this paper, we discuss the computation of the rank $i_k$ of homotopy group $\pi_k(G(A))$. For a large class of Kac-Moody groups, we construct Lie algebras with grade from the Poincar\'{e} series of their flag manifolds. And we interpret $i_{2k}$ as the dimension of the degree $2k$ homogeneous component of the Lie algebra we constructed. Since the computation of $i_{2k-1}$ is trivial, this gives a combinatorics interpretation of $i_k$ for all $k>0$.
\end{abstract}



\section{Introduction}

Let $A=(a_{ij})$ be an $n\times n$ integer matrix satisfying

\noindent (1) For each $i,a_{ii}=2$;

\noindent (2) For $i\not=j,a_{ij}\leq 0$;

\noindent (3) If $a_{ij}=0$, then $a_{ji}=0$,

\noindent then $A$ is called a Cartan matrix.

Let $h$ be the real vector space spanned by $\Pi^{\vee}=\{\alpha^{\vee}_1,\alpha^{\vee}_2,\cdots,\alpha^{\vee}_n\}$, denote the dual basis of $\Pi^\vee$ in the dual vector space $h^*$ by $\{\omega_1,\omega_2,\cdots,\omega_n\}$. That is $\omega_i(\alpha_j^{\vee})=\delta_{ij}$ for $1\leq i,j\leq n$. Let $\Pi=\{\alpha_1,\cdots,\alpha_n\}\subset h^*$ be given by $\langle\alpha^{\vee}_i,\alpha_j\rangle=a_{ij}$ for all $i,j$, then $\alpha_i=\sum\limits_{j=1}^n a_{ji}\omega_j$. The triple $(h,\Pi,\Pi^\vee)$ is called the realization of Cartan matrix $A$.
%Note that if the Cartan matrix $A$ is singular, then $\alpha_i,1\leq i\leq n$ is not a basis of $h^*$.
$\Pi$ and $\Pi^\vee$ are called respectively the simple root system and simple coroot system associated to Cartan matrix $A$.
%$\alpha_i, \alpha_i^\vee, \omega_i,1\leq i\leq n$ are called respectively the simple roots, simple coroots and fundamental dominant weights.

By the work of Kac\cite{Kac_68} and Moody\cite{Moody_68}, it is well known that for each Cartan matrix $A$, there is a Lie algebra $g(A)$ associated to
$A$ which is called the Kac-Moody Lie algebra.

%By the Proposition 1.6 of \cite{Kac_82}, the Cartan subgroup of $G(A)$ has rank $n$.

The Kac-Moody Lie algebra $g(A)$ is generated by $\alpha^\vee_i,e_i,f_i,1\leq i\leq n$ over $\mathbb{C}$,
with the defining relations:

(1) $[\alpha_i^{\vee},\alpha_j^{\vee}]=0$;

(2) $[e_i,f_j]=\delta_{ij}\alpha_i^{\vee}$;

(3) $[\alpha_i^{\vee},e_j]=a_{ij}e_j,[\alpha_i^{\vee},f_j]=-a_{ij}f_j$;

(4) $\mathrm{ad}(e_i)^{-a_{ij}+1}(e_j)=0, \mathrm{ad}(f_i)^{-a_{i j}+1}(f_j)=0$ for $i\not=j$.

For details, see Kac\cite{Kac_82} and Kumar\cite{Kumar_02}.

Kac and Peterson constructed the simply connected Kac-Moody group $G(A)$ with Lie algebra $g(A)$ in \cite{Kac_Peterson_83}\cite{Kac_Peterson_84}\cite{Kac_85}.

A Cartan matrix $A$ is symmetrizable if there exists an invertible diagonal matrix $D$ and a symmetric matrix $B$ such that $A=DB$. $g(A)$ and $G(A)$ are symmetrizable if the Cartan matrix $A$ is symmetrizable. A Cartan matrix $A$ is indecomposable if $A$ can not be decomposed into the sum $A_1\oplus A_2$ of two Cartan matrices $A_1,A_2$. The Kac-Moody Lie algebra or the Kac-Moody group is indecomposable if the Cartan matrix $A$ is indecomposable.

Indecomposable Cartan matrices and their associated Kac-Moody Lie algebras, Kac-Moody groups are divided into three types.

(1) Finite type, when $A$ is positive definite. In this case, $G(A)$ is just the simply connected complex semisimple Lie group with Cartan matrix $A$.

(2) Affine type, when $A$ is positive semi-definite and has rank $n-1$.

(3) Indefinite type otherwise.

For the Kac-Moody Lie algebra $g(A)$, there is Cartan decomposition $g(A)=h\oplus \sum\limits_{\alpha\in \Delta} g_{\alpha}$, where $h$ is a Cartan sub-algebra and $\Delta$ is the
root system. Let $b=h\oplus \sum\limits_{\alpha\in \Delta^+} g_{\alpha}$ be the Borel sub-algebra, then $b$ corresponds to a Borel subgroup $B(A)$ in the Kac-Moody group $G(A)$. The
homogeneous space $F(A)=G(A)/B(A)$ is called the flag manifold. By Kumar\cite{Kumar_02}, $F(A)$ is an ind-variety.

By rational homotopy theory\cite{Sullivan_77} the ranks of homotopy groups of $G(A)$ can be computed from the rational cohomology of $G(A)$. The rational cohomology rings of Kac-Moody groups and their flag manifolds of finite or affine type have been extensively studied by many people. For reference, see Pontrjagin\cite{Pontryagin_35}, Hopf\cite{Hopf_41}, Borel\cite{Borel_53_1}\cite{Borel_53}\cite{Borel_54}, Bott and Samelson\cite{Bott_Samelson_55}, Bott\cite{Bott_56}, Milnor and Moore\cite{Milnor_Moore_65} etc. The structure theory of cohomology rings is well established. For the indefinite case, there are some works by Kumar\cite{Kumar_85}, Kac\cite{Kac_85}, Kostant and Kumar\cite{Kostant_Kumar_86} and Kichiloo\cite{Kitchloo_98}, Zhao and Jin\cite{Zhao_Jin_12_4}. The fundamental structure and the explicit algorithm to determine the cohomology are founded. But except for the examples in \cite{Kitchloo_98}\cite{Zhao_Jin_12_4}, there are no concrete computational examples. This paper will work in this direction and give more examples.

%Below is the well known theorem about the rational cohomology of Hopf spaces.

%\noindent {\bf Theorem}(Hopf): Let $G$ be a connected H-space which has the homotopy type of a CW-complex, then the rational cohomology rings $H^*(G)$ is a Hopf algebra and as algebra it is %isomorphic to the the tensor product of a polynomial algebra $P(V_0)$ and
%a exterior algebra $\Lambda(V_1)$, where $V_0$ and $V_1$ are
%respectively the set of even and odd degree free generators of $H^*(G)$.

By Kichiloo\cite{Kitchloo_98} and Kumar\cite{Kumar_02}, it follows that the rational cohomology rings $H^*(G(A))$ and $H^*(F(A))$ are locally finite and generated by countable number of generators.

By the well known theorem of Hopf about the structure of the rational cohomology ring of a Hopf space $G$, we know $H^*(G)$ is a Hopf algebra and as algebra it is isomorphic to the tensor product of a polynomial algebra $P(V_0)$ and
an exterior algebra $\Lambda(V_1)$, where $V_0$ and $V_1$ are
respectively the set of even and odd degree free generators of $H^*(G)$.
Therefore the Poincar\'{e} series of the Kac-Moody group $G(A)$ of form
\begin{equation}\label{1}
 P_{G(A)}(q)=\prod\limits_{k=1}^{\infty}\frac{ (1-q^{2k-1})^{i_{2k-1}}}{(1-q^{2k})^{i_{2k}}}.
\end{equation}
\noindent By \cite{Sullivan_77} $i_k$ is the rank of homotopy group $\pi_k(G(A))$. The rational cohomology ring $H^*(G(A))$(even the rational homotopy type) is determined by the sequence $i_1,i_2,\cdots,i_k,\cdots$.

%\noindent{\bf Theorem 3: }For an indecomposable and indefinite Cartan matrix $A$, if $A$ is symmetrizable, then
%$$H^*(G(A))\cong \Lambda_{\mathbb{Q}}(y_3)\otimes \mathbb{Q}[z_1,\cdots,z_k,\cdots]$$ and $$H^*(F(A))\cong  \mathbb{Q}[\omega_1,\cdots,\omega_n]/<\psi>\otimes \mathbb{Q}[z_1,\cdots,z_k,\cdots]. %$$
%If $A$ is non symmetrizable, then
%$$H^*(G(A))\cong \mathbb{Q}[z_1,\cdots,z_k,\cdots]$$ and $$H^*(F(A))\cong  \mathbb{Q}[\omega_1,\cdots,\omega_n]\otimes \mathbb{Q}[z_1,\cdots,z_k,\cdots]. $$
%where $\deg z_k\geq 4$ are even for all $k$ and their degrees can be determined from the Poincar\'{e} series $P_A(q)$ and $\epsilon(A)$. By \cite{Zhao_Jin_12_3}, $y_3$ is corresponds to the %$W(A)$ invariant bilinear form $\psi$ on $h$.

Set $\epsilon(A)=1$ or $0$ depending on $A$ is symmetrizable or not as in Kac\cite{Kac_852}. By the results in \cite{Zhao_Jin_12_3}\cite{Zhao_Jin_12_4}, we have

\noindent{\bf Theorem 1: }The Poincar\'{e} series of $F(A)$ is
\begin{equation}\label{2}
       P_{F(A)}(q)=\frac{\prod\limits_{k=1}^{\infty} (1-q^{2k})^{i_{2k-1}}}{(1-q^2)^n}   \frac{1}{\prod\limits_{k=1}^{\infty}(1-q^{2k})^{i_{2k}}}.
\end{equation}

\noindent{\bf Theorem 2: }The sequence $i_1,i_2,\cdots,i_{k},\cdots$ can be computed from $P_{F(A)}(q)$ and $\epsilon(A)$. In particular $i_1=i_2=0,i_3=\epsilon(A)$ and $i_{2k-1}=0$ for $k\geq 3$.

For the computation of $i_{2k}$, see \cite{Jin_Zhao_11}\cite{Zhao_Jin_12_3}\cite{Zhao_Jin_12_4}.

%our main goal is the computation of sequence $i_1,i_2,\cdots,i_{k},\cdots$.
In this paper we will give a combinatorics realization of the rank $i_{2k}$ of homotopy group $\pi_k(G(A))$ for a large class of Kac-Moody groups. Since the computation for finite and affine cases has been obtained and for a decomposable Cartan matrix $A=A_1\oplus A_2$, $G(A)\cong G(A_1)\times G(A_2)$, we only consider the indecomposable and indefinite case.

%In fact $i_{2k}$ can be got by counting the dimension of degree $2k$ component of a certain finite generated Lie algebra $L(A)$(with grade) constructed from $A$. The Lie algebra $L(A)$ is %constructed from the Poincar\'{e} series $P_{F(A)}(q)$ for $G(A)$ with $rank(A)=3$ or $4$.

The content of this paper is arranged as follows. In section 2, we give some results about the Hilbert series of graded associative algebras which will be used in the later sections. In section 3 we construct a Lie algebra $L(A)$ with grade for a Kac-Moody group $G(A)$ with certain good property. And we interpret $i_{2k}$ as the dimension of the degree $2k$ component of $L(A)$. Since the realization is based on the Poincar\'{e} series of $F(A)$ we discuss the computation of Poincar\'{e} series of $F(A)$ in section 4. We are particularly interested in the case when the rank of Cartan matrix is $3$ or $4$. In section 5 we give some examples to show how our interpretation is implemented. In the last section we make a conjecture about the structure of grade Lie algebra $\pi_*(G(A))$ for certain type of Cartan matrix $A$.

\section{Hilbert series of graded associative algebras}
In this section we give some algebraic preparation. Our main reference for this section is Anick\cite{Anick_82}.
\subsection{Hilbert series of free product}
Let $A$ be a graded associative algebra over a field $K$, then $A=\sum\limits_{i=0}^\infty A_k$, where $A_k$ is the homogeneous degree $k$ component of $A$. Denote the augmented ideal $\sum\limits_{i=1}^\infty A_k$ of $A$ by $\widetilde A$. For two graded associative algebras $A_1,A_2$, denote by $A_1*A_2$ the free product of $A_1$ and $A_2$. In this paper we consider only connected associative algebra $A$. That is $A_0\cong K$. The Hilbert series of $A$ is $H_A(q)=\sum\limits_{i=0}^\infty q^k \dim A_k$. The Hilbert series satisfy the following property.

\noindent{\bf Lemma 1}(Lemaire\cite{Lemaire_74})Let $A_1,A_2$ be two connected graded associative algebras with Hilbert series $H_{A_1},H_{A_2}$, then the Hilbert series of the free product $A=A_1*A_2$ satisfies $\displaystyle{\frac{1}{H_{A}}=\frac{1}{H_{A_1}}+\frac{1}{H_{A_2}}-1}$.

\noindent {\bf Proof: }Let $S_{1},S_2$ be the set of additive basis of augmented ideal $\widetilde A_1,\widetilde A_2$ respectively. Since $A$ is the free product of $A_1,A_2$, we can construct a canonical basis $S$ of $\widetilde A$ from $S_1,S_2$ whose elements are of the form of finite product $\alpha_{i_1}\alpha_{i_2}\cdots\alpha_{i_k},k>0$ such that if $\alpha_{i_j}\in S_1$ then $\alpha_{i_{j+1}}\in S_2$ and if $\alpha_{i_j}\in S_2$ then $\alpha_{i_{j+1}}\in S_1$ for $1\leq j\leq k-1$. Let $F$ be the subspace of $\widetilde A$ spanned by those elements of $S$ starting from $\alpha_{i_1}\in S_1$ and $G$ be the subspace of $\widetilde A$ spanned by those elements of $S$ starting from $a_{i_1}\in S_2$. Then $A\cong K\oplus F \oplus G$.

Considering the Hilbert series of the two sides, we get
$$H_A=1+H_F+H_G.$$

It is obvious that $F=\widetilde A_1 A$ and $G=\widetilde A_2 A$. Hence $A=K\oplus \widetilde{A}_1 A\oplus \widetilde{A}_2 A$. We called this formula the first order expansion of $A$. It is easy to check that we have the following second order expansion of $A$
$$A=K\oplus \widetilde{A}_1 \oplus \widetilde{A}_2 \oplus \widetilde{A}_1 \widetilde{A}_2\oplus \widetilde{A}_2 \widetilde{A}_1 \oplus \widetilde{A}_1 \widetilde{A}_2 F \oplus \widetilde{A}_2 \widetilde{A}_1 G.$$
So
$$H_A=1+H_{A_1}-1+H_{A_2}-1+2(H_{A_1}-1)(H_{A_2}-1)+(H_{A_1}-1)(H_{A_2}-1)(H_F+H_G-1)$$
Simplifying this formula, we get $$\displaystyle{\frac{1}{H_{A}}=\frac{1}{H_{A_1}}+\frac{1}{H_{A_2}}-1}.$$

Let $T(x_1,\cdots,x_m)$ be the tensor algebra generated by $x_1,\cdots ,x_m$. Since $T(x_1,\cdots,x_m)\cong T(x_1)*\cdots * T(x_m)$ and for $x$ with $\deg x=d, \displaystyle{H_{T(x)}=\frac{1}{1-q^d}}$. We have

\noindent{\bf Corollary 1: }For tensor algebra $A=T(x_1,\cdots,x_m)$ with $\deg x_i=d_i,1\leq i\leq m$, then $$\displaystyle{H_A=\frac{1}{1-q^{d_1}-\cdots-q^{d_m}}}.$$

\subsection{Strongly free set}
Let $A$ be a graded associative algebra and $B$ be a subalgebra of $A$, then the quotient homomorphism $\pi:A\to A/ABA$ is surjective. Let $\rho: A/ABA\to A$ be a chosen linear section of $\pi$, then there is a homomorphism $\mathrm{id} * \rho:B * (A/ABA)\to A$.

The following definition of strongly free set can be regarded as the generalization of the concept of regular sequences for commutative algebras to non-commutative algebras.

\noindent{\bf Definition 1}(Anick\cite{Anick_82})Let $A$ be a graded associative algebra and $B$ be a subalgebra of $A$, $B$ is called a weak summand of $A$ if the homomorphism $\mathrm{id} * \rho:B * (A/ABA)\to A$ is an isomorphism of $K$-vector spaces. Let $\alpha=\{\alpha_1,\cdots,\alpha_k\}$ be a graded set in $A$, $\alpha$ is called a strongly free set in $A$ if the subalgebra $K\langle \alpha\rangle $ generated by $\alpha$ in $A$ is a free algebra and $K\langle \alpha\rangle$ is a weak summand of $A$.

Let $\alpha$ be a strongly free set in $A$, $\alpha$ generates an ideal $A\alpha A\subset A$ in $A$, the following lemma gives the relation between the Hilbert series $H_A$ and $H_{A/A\alpha A}$.

\noindent{\bf Lemma 2: }Let $A$ be a connected graded associative algebra and $\alpha=\{\alpha_1,\cdots,\alpha_n\}$ be a strongly free set in $A$. If the degrees of elements in $\alpha$ are $e_1,\cdots,e_n$, then $\displaystyle{\frac{1}{H_{A/A\alpha A}}=\frac{1}{H_A}+q^{e_1}+\cdots+q^{e_n}}$.

\noindent{\bf Proof: }Since $\mathrm{id} * \rho:T(\alpha) * (A/ABA)\to A$ is an isomorphism of vector spaces, we get $$\displaystyle{\frac{1}{H_{A}}=\frac{1}{H_{A/A\alpha A}}+\frac{1}{H_{K\langle \alpha\rangle}}-1}$$
Combing with $\displaystyle{H_{K\langle \alpha\rangle }=\frac{1}{1-e^{e_1}-\cdots-q^{e_l}}}$, we prove the lemma.

%\noindent{\bf Remark: }The result is also true for strongly free set with infinite cardinal.

For a connected graded associative algebra $A$ generated by set $X=\{x_1,\cdots,x_m\}$, then $$A\cong T(x_1,\cdots,x_m)/I$$ where $I$ is the ideal of relation of $A$ with respect to generators set $X$. If $I$ is generated by a strongly free set in $T(x_1,\cdots,x_m)$, then we have the following result.

\noindent{\bf Corollary 2: }Let $A$ be a connected graded associative algebra with generator set $X=\{x_1,\cdots,x_m\}$ and relation set $\alpha=\{\alpha_1,\alpha_2,\cdots,\alpha_n\}$. If $\alpha$ is a strongly free set, then the Hilbert series of $A$ is $$H_A=\displaystyle{\frac{1}{1-q^{d_1}-\cdots-q^{d_m}+q^{e_1}+\cdots+q^{e_n}}}.$$
\noindent where $d_1,d_2\cdots,d_m$ are the degrees of set $\{x_1,\cdots,x_m\}$ and $e_1,e_2\cdots,e_n$ are the degrees of set $\{\alpha_1,\alpha_2,\cdots,\alpha_n\}$.

\subsection{Hilbert series of universal enveloping algebra of Lie algebra with grade}
Let $L=\bigoplus \limits_{i=1}^\infty L_i$ be a Lie algebra and $[L_k,L_l]\subset L_{k+l}$ for all $k,l>0$, then we say $L$ is a Lie algebra with grade. For a Lie algebra $L$ with grade, its universal enveloping algebra $U(L)$ is an graded associative Hopf algebra. The coproduct on $U(L)$ is defined by $\delta(x)=1\otimes x+x\otimes 1$ for $x\in L$ and $\delta$ is cocommutative. Therefore By \cite{Milnor_Moore_65} $U(L)$ is primitively generated. For a free Lie algebra $L$ with grade generated by $X=\{x_1,x_2,\cdots,x_m\}$, the universal enveloping algebra is the tensor algebra $T(x_1,x_2,\cdots,x_m)$.

\noindent{\bf Lemma 3: }Let $L=\bigoplus\limits_{k=1}^\infty L_k$ be a Lie algebra with grade and $j_k=\dim L_k$,
%Suppose $U(L)$ be the universal enveloping algebra of $L$,
then the Hilbert series $H_{U(L)}=\displaystyle{\frac{1}{\prod\limits_{k=1}^{\infty}{(1-q^{k})^{j_{k}}}}}$.

This lemma is derived directly from the Poincar\'{e}-Birkhoff-Witt Theorem.
%\noindent{\bf Proof: }
%Since the coproduct on $U(L)$ is cocommutative and primitively generated. The dual Hopf algebra $U(L)^*$ is commutative. In fact it is a polynomials algebra with grade., we have $H_{U(L)^*}=\displaystyle{\frac{1}{\prod\limits_{k=1}^{\infty}{(1-q^{k})^{i_{k}}}}}$. Since we have
%$H_{U(L)}=H_{U(L)^*} $,

Let $L$ be a Lie algebra with grade and $\alpha\subset A$ be a graded set and $J$ be the quotient Lie algebra of $L$ with respect to the ideal $I$ generated by $\alpha$, then we have

\noindent{\bf Lemma 4: }The universal enveloping algebra $U(J)$ is isomorphic to the quotient Hopf algebra of $U(L)$ with respect to the ideal $U(L)IU(L)$.

For the proof of this lemma, see Bourbaki \cite{Bourbaki}.
%\noindent{\bf Lemma : }Let $L$ be a graded free Lie algebra generated by graded set $X=\{x_1,\cdots,x_l\}$ with degree of $x_i$ be $d_i$ and $i_k=\dim L_k$,  then $\frac{1}{1-q^{d_1}-\cdots-q^{d_l}}=\frac{\prod\limits_{k=1}^{\infty}{(1-q^{2k-1})^{i_{2k-1}}}}{\prod\limits_{k=1}^{\infty}{(1-q^{2k})^{i_{2k}}}}$.

\noindent{\bf Definition 2: }Let $L$ be a Lie algebra with grade, $\alpha=\{\alpha_1,\alpha_2,\cdots,\alpha_m\}\subset L$ is called a strongly free set in $L$ if and only if the image of $\alpha$ in $U(L)$ is a strongly free set.

\noindent{\bf Lemma 5: }Let $L$ be a Lie algebra with grade generated by graded set $X=\{x_1,\cdots,x_m\}$ with defining relation set $\alpha=\{\alpha_1,\cdots,\alpha_n\}$.
If $\alpha$ is strongly free set, then the Hilbert series of $U(L)$ is $$\displaystyle{\frac{1}{1-q^{d_1}-\cdots-q^{d_m}+q^{e_1}+\cdots+q^{e_n}}=\frac{1}{\prod\limits_{k=1}^{\infty}{(1-q^{k})^{j_{k}}}}}.$$
\noindent where $d_1,d_2\cdots,d_m$ are the degrees of set $\{x_1,\cdots,x_m\}$ and $e_1,e_2\cdots,e_n$ are the degrees of set $\{\alpha_1,\alpha_2,\cdots,\alpha_n\}$.

\section{Rank of homotopy group $\pi_k(G(A))$}
\subsection{The theorem of Milnor and Moore}
In this section we use the results in previous section to discuss the rank $i_k$ of homotopy group $\pi_k(G(A))$.

For a connected Hopf space $G$ with unit and homotopy associative multiplication, the homology $H_*(G)$ is a Hopf algebra.
The diagonal $\Delta: G \to G\times G$ is co-commutative,
implying that the coproduct of $H_*(G)$ is commutative. Hence by \cite{Milnor_Moore_65} $H_*(G)$ is primitively generated.

On the rational homotopy group $\pi_*(G)$, the Samelson product $[\  , \ ]:\pi_p(G)\times \pi_q(G)\to \pi_{p+q}(G)$ is defined as
$$[\alpha,\beta](s\wedge t)=\alpha(s)\beta(t)\alpha(s)^{-1}\beta(t)^{-1},s\in S^p,t\in S^q.$$
$\pi_*(G)$ forms a graded Lie algebra over $\mathbb{Q}$ with Samelson product.

\noindent{\bf Theorem}(Milnor-Moore) Let $G$ be a connected homotopy associative H-space with unit and $\chi: \pi_*(G)\to H_*(G)$ be the Hurewicz morphism of graded Lie algebras, then the induced morphism $\widetilde \chi: U(\pi_*(G))\to H_*(G)$ is an isomorphism of Hopf algebras. Where $U(\pi_*(G)))$ is the enveloping algebra of $\pi_*(G)$.

By this theorem, to determine the rank $i_k$ of $\pi_k(G(A))$, we only need to consider the Hilbert series of $H_*(G(A))$. Since the Hilbert series of $H_*(G(A))$ contains the same information as the Poincar\'{e} series of $G(A)$, we discuss the Poincar\'{e} series of $G(A)$ for convenience. $H_*(G(A))$ is the tensor product of a polynomial algebra with even degree generators and an exterior algebra with odd degree generators. In this paper we only consider indecomposable and indefinite Cartan matrix $A$. In this case it is proved in \cite{Zhao_Jin_12_3} that if $A$ is symmetrizable then the exterior algebra part of $H^*(G(A))$ is generated by one degree $3$ generators and if $A$ is not symmetrizable then $H^*(G(A))$ has no exterior algebra part.

\subsection{%Poincar\'{e} series of
Chow ring of $G(A)$}
Lie algebra $\pi_*(G(A))$ is a graded Lie algebra whose universal enveloping algebra is $H_*(G(A))$. For the Lie sub-algebra $\pi_{even}(G(A))=\sum\limits_{i=1}^\infty \pi_{2k}(G(A))$, the universal enveloping algebra is $H_{even}(G(A))$. The dual Hopf algebra of $H^{even}(G(A)$ is isomorphism to the Chow ring $\mathrm{Ch}^*(G(A))$. As algebra $\mathrm{Ch}^*(G(A))$ is the subalgebra of $H^*(G(A))$ generated by even dimensional generators of degree great than $2$. By relating $\pi_{even}(G(A))$ with Chow ring of $G(A)$ we transform the computation of $i_{2k}$ to the computation of Hilbert series of Chow ring.

%By this way we needn't to consider graded Lie algebra $\pi_*(G(A))$ and use Lie algebra $\pi_{even}(G(A))$, which can be simplified as a Lie algebra with grade. This way we transform the problem to consider the graded Lie algebras to Lie algebras with grade.
%We avoid to introduce the concept of

\noindent{\bf Lemma 6: }The Hilbert series of $\mathrm{Ch}^*(G(A))$ is $$C_A(q)=P_{F(A)}(q)(1-q^2)^n(1-q^4)^{-\epsilon(A)}$$ and $$\hspace*{5.0cm} \displaystyle{C_A(q)=\prod\limits_{k=2}^{\infty}\frac{ 1}{(1-q^{2k})^{i_{2k}}}}    \hspace*{5.0cm}   (3)$$
For reference see Kac\cite{Kac_852}.

%\subsection{Strongly positive polynomials}

%Suppose polynomial
%$Q(q)=1-a_1 q^{ d_1}-\cdots-a_m q^{d_m}+b_1 q^{ e_1}+\cdots+b_n q^{e_n}$ with $a_i>0,1\leq i\leq m; b_j>0,1\leq j\leq n$ and $d_1<d_2<\cdots<d_m,e_1<e_2<\cdots<e_n$.

%Let the free graded associative algebra with generators of degree $d_1,d_2,\cdots,d_m$ and the number of degree

%\noindent{\bf Definition 3: }The polynomial $Q(q)$ is called a strongly positive polynomial if $\displaystyle{\frac{1}{Q(q)}}$ is Hilbert series of the quotient algebra of a free graded associative algebra $A$ with respect to the ideal generated by a strongly free set $\alpha\in A$, such that the generators of $A$ have degree $d_1,d_2,\cdots,d_m$ and the number of degree $d_i$ generators is $a_i$ for all $i$, the degrees of elements in $\alpha$ are $e_1,e_2,\cdots,e_n$ and the number of degree $d_i$ elements is $b_j$ for all $j$.

%Denoted by $Q(q)\succ 0$. For two polynomial $f(q),g(q)$ with constant item $0$, we say $f(q)\prec g(q)$ if $\displaystyle{\frac{1}{1-f(q)+g(q)}}$ is strongly positive polynomial.

%\noindent {\bf Lemma 7: } If $Q_1,Q_2$ are strongly positive, then $Q_1+ Q_2-1$ is strongly positive.

%The lemma is proved by using Lemma 2.

%the observation $result The relation $\prec$ on the set of all polynomials with nonnegative coefficients and constant item $0$ is a partial order relation.

%If $Q(q)=1-f(q)+g(q)$ is strongly positive, then the expansion coefficients of power series $\displaystyle{\frac{1}{1-f(q)+g(q)}}$ is nonnegative.

%We need the following more stronger definition.

%

\subsection{Realization of $i_k$}

The Poincar\'{e} series $P_{F(A)}(q)$ is of the form $\displaystyle{\frac{\prod\limits_{i=1}^r [t_i]}{Q(q^2)}}$, where $Q(q)$ is a polynomial of $q$ with constant item $1$ and $[d_i]=\displaystyle{\frac{1-q^{2t_i}}{1-q^2}}$. By \cite{Jin_Zhao_11} $\prod\limits_{i=1}^r [t_i]$ is in fact the least common multiple of those Poincar\'{e} series of flag manifolds associated to the finite type principal sub-matrices of $A$. We assume the polynomial
$Q(q)$ is of the form $1-a_1 q^{ d_1}-\cdots-a_m q^{d_m}+b_1 q^{ e_1}+\cdots+b_n q^{e_n}$ with $a_i>0,1\leq i\leq m; b_j>0,1\leq j\leq n$ and $d_1<d_2<\cdots<d_m,e_1<e_2<\cdots<e_n$.

We have $$\displaystyle{C_A(q)=\prod\limits_{k=2}^{\infty}\frac{ 1}{(1-q^{2k})^{i_{2k}}}}.$$
and by Zhao-Jin\cite{Zhao_Jin_12_4} there exists a unique sequence $j_1,j_2,\cdots,j_k\cdots$ such that 
$$\displaystyle{\frac{1}{Q(q)}=\prod\limits_{k=1}^{\infty}\frac{ 1}{(1-q^{k})^{j_{k}}}}. $$
Substituting the above two formulas into 
$$\hspace*{2cm} C_A(q)=\displaystyle{\frac{(1-q^2)^{n-r}(1-q^4)^{-\epsilon(A)}\prod\limits_{i=1}^r (1-q^{2t_i})}{Q(q^2)}}.\hspace*{3cm}  (4)$$
We get:
%And by We called the sequence characteristic sequence of $\displaystyle{\frac{1}{Q(q)}}$. Then we have
%$$\displaystyle{\frac{1}{Q(q^2)}=\prod\limits_{k=1}^{\infty}\frac{ 1}{(1-q^{2k})^{j_{k}}}}.$$
%And

\noindent{\bf Lemma 7: }The series $i_k$ satisfy $i_2=0,i_4=j_2-l_2+\epsilon(A),i_{2k}=j_k-l_k,k>2$, where $l_k=\#\{i|t_i=k,1\leq i \leq r\}$.

We give the following definition.

\noindent{\bf Definition 3: }A polynomial $Q(q)$ is called a strongly positive polynomial if there exists a free Lie algebra $L$ with grade and a strongly free set $\alpha\in L$, such that $\displaystyle{\frac{1}{Q(q)}}$ is the Hilbert series of the quotient algebra of the universal enveloping algebra $U(L)$ with respect to the ideal generated by $\alpha\in L$. A Kac-Moody groups is called a good Kac-Moody group if it correspond to a strongly positive polynomials $Q(q)$¡£  

A large class of Kac-Moody groups are good. %Assume the generators of $L$ have degree $d_1,d_2,\cdots,d_m$ and the number of degree $d_i$ generators is $a_i$ for all $i$, the degrees of elements in $\alpha$ are $e_1,e_2,\cdots,e_n$ and the number of degree $d_i$ elements is $b_j$ for all $j$.
For a good Kac-Moody groups $G(A)$, there exists a Lie algebra $L(A)$ such that the Hilbert series of $U(L)$ is $Q(q)$. In this case $j_k$ is just the dimension of degree $k$ homogeneous component of $L(A)$. So we have

\noindent{\bf  Theorem 3: }For a good Kac-Moody groups $G(A)$, the rank $i_k$ of the homotopy group $\pi_k(G(A))$ satisfy $i_2=0,i_4=j_2-l_2+\epsilon(A),i_{2k}=j_k-l_k,k>2$, where $j_k$ is the rank of the degree $k$ component of Lie algebra $L(A)$ and $l_k=\#\{i|t_i=k,1\leq i \leq r\}$.

\section{The computation of Poincar\'{e} series $P_{F(A)}(q)$}
We need an algorithm to compute the Poincar\'{e} series of flag manifolds.
\subsection{General results about the Poincar\'{e} series of flag manifolds}

The Weyl group $W(A)$ associated to a Cartan matrix $A$ is the group generated by the Weyl reflections $\sigma_i:h^*\to h^*$ with respect to simple co-roots $\alpha_i^{\vee},1\leq i\leq n$, where $\sigma_i(\alpha)=\alpha-\langle \alpha,\alpha_i^\vee\rangle\alpha_i$. $W(A)$ has a Coxeter presentation
$$W(A)=<\sigma_1,\cdots,\sigma_n|\sigma^2_i=e,1\leq  i\leq n; (\sigma_i\sigma_j)^{m_{ij}}=e,1\leq i<j\leq n>.$$ where $m_{ij}=2,3,4,6 $ or $\infty$ as $a_{ij}a_{ji}=0,1,2,3 $ or $\geq 4$ respectively. For details see Kac\cite{Kac_82}, Humphreys\cite{Humphreys_90}.

%The action of $\sigma_i$ on fundamental dominant weights is given by $\sigma_i(\omega_j)=\omega_j-(\omega_j,\alpha_i^{\vee})\alpha_i=\omega_j-\delta_{ji}\alpha_i$.

Each element $w\in W(A)$ has a decomposition of the form $w=\sigma_{i_1}\cdots \sigma_{i_k},1\leq i_1,\cdots,i_k\leq n$. The length of $w$ is defined as the least integer $k$ in all of those decompositions of $w$, denoted by $l(w)$. The Poincar\'{e} series of $g(A)$ is the power series $P_A(q)=\sum\limits_{w\in W(A)} q^{2l(w)}$. 

%Hence $P_A(q)$ only depends on the structure of Weyl group $W(A)$ and the length function on it.

%The flag manifold $F(A)$ admits a CW-decomposition of Schubert cells which are indexed by the elements of Weyl group $W(A)$. For each $w\in W(A)$, the real dimension of Schubert variety $X_w$ is $2l(w)$. So the Poincar\'{e} series of flag manifold $F(A)$ is just the Poincar\'{e} series $P_A(q)$ of $g(A)$.

Steinberg\cite{Steinberg_68} proved that
the Poincar\'{e} series $P_A(q)$ of the flag manifold $F(A)$ of a Lie group $G(A)$ is a rational function. And the result is easy to extend to the Poincar\'{e} series of a general Kac-Moody group $G(A)$ or its flag manifolds $F(A)$.

In \cite{Jin_Zhao_11} the authors discussed the computation of Poincar\'{e} series of flag manifolds of Kac-Moody flag groups. Let $A$ be an $n\times n$ Cartan matrix, $S=\{1,2,\cdots,n\}$. For each $I\subset S$, let $A_I$ be the Cartan matrix $(a_{ij})_{{i,j}\in I}$. Let $\dim A_I$ be the complex dimension of the flag manifolds $F(A_I)$. The following lemma is used in this paper.

\noindent{\bf Lemma 8}(Steinberg\cite{Steinberg_68})Let $A$ be an indefinite Cartan matrix, then we have
$$ \hspace*{5cm}\sum\limits_{I\subset S}(-1)^{|I|}\dfrac{P_{F(A)}(q)}{P_{F(A_I)}(q)}=0 \hspace*{5.7cm}(5)$$ and $$\hspace*{4cm}\displaystyle{P_{F(A)}(q^{-1})=\sum\limits_{I\subset S,\dim A_I < \infty}\frac{(-1)^{|I|}}{P_{F(A_I)}(q)}}\hspace*{4cm} (6)$$

The Poincar\'{e} series $P_{F(A)}(q)$ can be computed through Steinberg's formula by a recurrence procedure.

\subsection{Poincar\'{e} series of flag manifolds of rank 3 and 4}
By the results in \cite{Jin_10}\cite{Jin_Zhao_11}, for a Cartan matrix $A$ the Poincar\'{e} series $P_{F(A)}(q)$ is determined by the Coxeter graph $\Gamma(A)$. For a Coxeter graph $\Gamma(A)$ we define the reduced graph is the graph obtained by replacing all the $k$-fold edges between pairs of vertices by one-fold edge.

\noindent{\bf Lemma 9: }Let $A$ be an indecomposable rank $3$ Cartan matrix, then its reduced Coxeter graph is of the following two types.

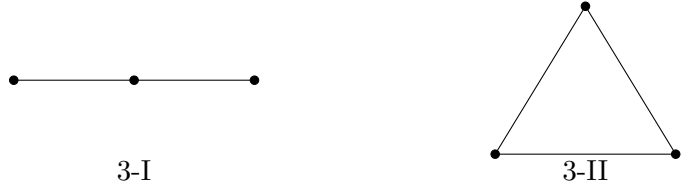
\begin{center}
\begin{tikzpicture}[scale=0.8]
{\filldraw [black] (10,0.5) circle (2pt);}
{\filldraw [black] (12,0.5) circle (2pt);}
{\filldraw [black] (14,0.5) circle (2pt);}

{\draw (10,0.5) -- +(2,0);}
{\draw (12,0.5) -- +(2,0);}
\draw (12,-1)  node{3-I};

{\filldraw [black] (18,-0.73) circle (2pt);}
{\filldraw [black] (21,-0.73) circle (2pt);}
{\filldraw [black] (19.5,1.73) circle (2pt);}

{\draw (18,-0.73) -- +(1.5, 2.46);}
{\draw (21,-0.73) -- +(-1.5, 2.46);}
{\draw (18,-0.73) -- +(3, 0);}
\draw (19.5,-1)  node{3-II};
\end{tikzpicture}
\end{center}

\noindent{\bf Lemma 10: }Let $A$ be an indecomposable rank $4$ Cartan matrix, then its Coxeter graph is of the following five types.
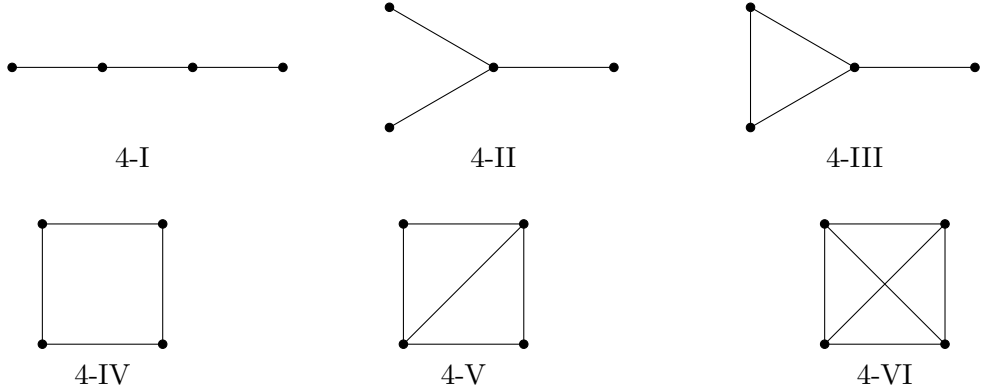
\begin{center}
\begin{tikzpicture}[scale=0.8]
\foreach \x in {10,11.5,13,14.5}
\foreach \y in {0.5}
{\filldraw [black] (\x,\y) circle (2pt);}

\foreach \x in {10,11.5,13}
\foreach \y in {0.5}
{\draw (\x,\y) -- +(1.5,0);}

\draw (12,-1)  node{4-I};

{\filldraw [black] (10+8,0.5) circle (2pt);}
{\filldraw [black] (12+8,0.5) circle (2pt);}
{\filldraw [black] (10+8-1.73,-0.5) circle (2pt);}
{\filldraw [black] (10+8-1.73,1.5) circle (2pt);}

{\draw (10+8,0.5) -- +(-1.73, 1);}
{\draw (10+8,0.5) -- +(-1.73, -1);}
{\draw (10+8,0.5) -- +(2,0);}
\draw (18,-1)  node{4-II};

{\filldraw [black] (10+14,0.5) circle (2pt);}
{\filldraw [black] (12+14,0.5) circle (2pt);}
{\filldraw [black] (10+14-1.73,-0.5) circle (2pt);}
{\filldraw [black] (10+14-1.73,1.5) circle (2pt);}

{\draw (10+14,0.5) -- +(-1.73, 1);}
{\draw (10+14,0.5) -- +(-1.73, -1);}
{\draw (10+14,0.5) -- +(2,0);}
{\draw (10+14-1.73,-0.5) --(10+14-1.73,1.5);}
\draw (24,-1)  node{4-III};

\end{tikzpicture}
\end{center}

\begin{center}
\begin{tikzpicture}[scale=0.8]

{\filldraw [black] (10+6,-0.5) circle (2pt);}
{\filldraw [black] (12+6,-0.5) circle (2pt);}
{\filldraw [black] (10+6,1.5) circle (2pt);}
{\filldraw [black] (12+6,1.5) circle (2pt);}

{\draw (10+6,-0.5) -- +(2, 0);}
{\draw (10+6,-0.5) -- +(0, 2);}
{\draw (12+6,1.5) -- +(-2, 0);}
{\draw (12+6,1.5) -- +(0, -2);}
\draw (17,-1)  node{4-IV};

{\filldraw [black] (10+8+4,-0.5) circle (2pt);}
{\filldraw [black] (12+8+4,-0.5) circle (2pt);}
{\filldraw [black] (10+8+4,1.5) circle (2pt);}
{\filldraw [black] (12+8+4,1.5) circle (2pt);}

{\draw (10+8+4,-0.5) -- +(2, 0);}
{\draw (10+8+4,-0.5) -- +(0, 2);}
{\draw (12+8+4,1.5) -- +(-2, 0);}
{\draw (12+8+4,1.5) -- +(0, -2);}
{\draw (18+4,-0.5) -- +(2, 2);}
\draw (23,-1)  node{4-V};

{\filldraw [black] (10+19,-0.5) circle (2pt);}
{\filldraw [black] (12+15+4,-0.5) circle (2pt);}
{\filldraw [black] (10+15+4,1.5) circle (2pt);}
{\filldraw [black] (12+15+4,1.5) circle (2pt);}

{\draw (10+15+4,-0.5) -- +(2, 0);}
{\draw (10+15+4,-0.5) -- +(0, 2);}
{\draw (12+15+4,1.5) -- +(-2, 0);}
{\draw (12+15+4,1.5) -- +(0, -2);}
{\draw (10+15+4,-0.5) -- +(2, 2);}
{\draw (10+15+6,-0.5) -- +(-2, 2);}
\draw (30,-1)  node{4-VI};
\end{tikzpicture}
\end{center}

The automorphisms of these graphs are: 3-I, $\mathbb{Z}_2$; 3-II, $S_3$; 4-I, $\mathbb{Z}_2$; 4-II, $S_3$; 4-III, $\mathbb{Z}_2$; 4-IV, ${D}_4$; 4-V, $\mathbb{Z}_2\times \mathbb{Z}_2$; 4-VI, $S_4$.

The computation results for the Poincar\'{e} series of $F(A)$ for rank $3$ and $4$ Cartan matrices are listed in the appendix A.
%\noindent{\bf Lemma : }For rank $3$ Cartan matrix of type I,

%\noindent{\bf Lemma : }For rank $3$ Cartan matrix of type I,

%The Poincar\'{e} series of $P_A(q)$ have the form of $\displaystyle{P_A(q)=\frac{\prod\limits_{i=1}^m [d_m]}{Q(q^2)}}$, where $Q(q)$ is in $\mathbb{Z}_1[q]$.

\section{Some examples}
\noindent{\bf Example 1: }For Cartan matrix $A=\left(
                                                 \begin{array}{ccc}
                                                   2 & -1 & -1 \\
                                                   -3 & 2 & -1 \\
                                                   -1 & -1 & 2 \\
                                                 \end{array}
                                               \right)$, $\epsilon(A)=0$. The Poincar\'{e} series of $F(A)$ is
$$P_{F(A)}(q)=\frac{[2][6]}{t^{12}-t^{10}-t^8+t^6-t^4-t^2+1}.$$

In this example $Q(q)=t^6-t^5-t^4+t^3-t^2-t+1$, we define a Lie algebra $L(A)=\langle x_1,x_2,x_4,x_5|[x_1,x_2],[x_1,x_5]\rangle$ with degree $i$ generator $x_i$ for $i=1,2,4,5$, then $[x_1,x_2],[x_1,x_5]$ form a strongly free set in $L$. Let the dimension of degree $k$ component of $L(A)$ be $j_k$, then by Theorem 3 we have $i_2=0,i_{4}=j_2-1,i_{12}=j_6-1$ and $i_{2k}=j_k$ for $k\not=1,2,6$.

\noindent{\bf Example 2: }For Cartan matrix $A=\left(
                                                 \begin{array}{cccc}
                                                   2 & -1 & -1 & -1 \\
                                                   -2 & 2 & -1 & -1  \\
                                                   -1 & -1 & 2 & -1  \\
                                                   -1 & -1 & -1 & 2
                                                 \end{array}
                                               \right)$, the $\epsilon(A)=0$. The Poincar\'{e} series of $F(A)$ is
$$P_{F(A)}(q)=\frac{[2][3][4]}{3t^{12}+t^{10}-t^8-t^6-3t^4-t^2+1}.$$

So $Q(q)=3t^6+t^5-t^4-t^3-3t^2-t+1$, we define a Lie algebra 

$L(A)=\langle x_1,x_{21},x_{22},x_{23},x_4|[x_1,x_4],[x_{21},x_4],[x_{22},x_4],[x_{23},x_4]\rangle$. 

Where $\deg x_1=1,\deg x_4=4,\deg x_{2i}=2, \forall i$. Then $[x_1,x_4],[x_{21},x_4],[x_{22},x_4],[x_{23},x_4]$ form a strongly free set in $L$ and we have $i_2=0,i_{4}=j_2-1,i_{6}=j_3-1, i_8=j_4-1$ and $i_{2k}=j_k$ for $k\not=1,2,3,4$.

%\section{The rational cohomology ring of $G(A)$}

\section{Criteria for strongly free set }

There is no easily applied criterion to determine whether or not
a given graded set $\alpha$ in an algebra $A$ is strongly free. But for free algebras  Anick gave some criteria in \cite{Anick_82}. We cite the corresponding results.

\noindent{\bf Definition 4: }Let $S$ be any locally finite graded set and let $B$ be the free
monoid on S. A set of monomials $\alpha=\{\alpha_1,\alpha_2,\cdots,\alpha_n \} \subset B- \{ 1\}$ is
combinatorially free iff (a) no $\alpha_i$ is a sub-monomial of $\alpha_j$ for $i\not= j$ and
(b) whenever $\alpha_i=x_1 y_1$, and $\alpha_j=x_2 y_2$ for $x_1, y_1, x_2, y_2\in B - \{ 1\}$ we have
$y_1\not=x_2$.

Condition (b) says that the beginning of one monomial cannot
be the same as the ending of another (or the same) monomial.

\noindent{\bf Theorem 4: }(Anick\cite{Anick_82}) Let $A = K(S)$, let $B$ be the free monoid on $S$ and suppose
$\alpha = \{\alpha_1,\alpha_2,\cdots,\alpha_n \}$ $\subset B- \{ 1\}$ is a set of monomials. Then $\alpha$ is strongly free in $A$
if and only if $\alpha$ is combinatorially free.

%Fix a graded set $S$ and let $B$ be the free monoid on $S$. For $x, y \in B$,
%we say that $y$ is a submonomial of $x$ iff $y = 1$ or $x =\alpha_1\alpha_2\alpha_t$ and $y =
%\alpha_{i_k}\alpha_{i_{k+1}}\cdots\alpha_{i_l}$ for some $1 \leq m \leq n \leq t$ and $\alpha_{i_j}\in S$. ¡°Submonomial of¡± is a partial
%ordering on B. An order ideal of monomials is a non-empty subset $M\subset B$
%which is closed under taking submonomials, i.e., if $x\in M$ and $y$ is a
%submonomial of $x$, then $y\in M$. Since $M$ is non-empty, $1 \in M$.
%The following lemma is virtually identical to the corresponding fact about
%commutative graded algebras.

%\noindent{\bf Lemma :}Let $S$ generate $A$ and let $f: k(S)\to A$ be a subjection.
%Let $B$ be the free monoid on $S$. Then there is an order ideal of monomials
%$M \subset B$ such that the elements $f (x), x\in M$, form a basis for $A$ as a $K$-module.

In the following, we give a monomials basis $M$ for a free monoid $B=K(S)$. If $S$ is empty, then $B= K$ and the only monomial is $1$. Otherwise choose any total
ordering on $S$. Define an ordering $e$ on $B$ as follows. For monomial$x, y\in B, x < y$ iff
$e(x)< e(y)$; if $e(x) = e(y)$, compare $x$ and $y$ using the lexicographic ordering
induced on $B$ by the ordering for $S$. Since $(S, e)$ is locally finite, $e^{-1}(n) \cap B$
is finite for each $n$, and $B$ is isomorphic as an ordered set to the positive
integers. This ordering has the additional property that if $u, w, x, y \in B$ and
$x < y$, then $uxw < vyw$.

Given a nonzero element $x\in K(S)$,
write $x$ as a linear combination of monomials $x = c_1 y_1 +\cdots + c_l y_l$, where
$c_i \in K$. If $y_i$ is the largest monomial(in the sense of the ordering of the monomials) for which $c_i\not=0$, then $y_i$ is called the high term of $x$.

\noindent{\bf Theorem 5: }(Anick\cite{Anick_82})Let $A= K(S)$ and suppose $\alpha = \{\alpha_1,\alpha_2,\cdots,\alpha_n \} \subset \widetilde A -\{ 0\}$.
Using any fixed ordering on $S$, let $\hat\alpha_i$ be the high term of $\alpha_i$ for each $\alpha_i\in \alpha$, then $\alpha$ is strongly free in $A$ if $\hat\alpha = \{\hat\alpha_1,\hat\alpha_2,\cdots,\hat\alpha_n \}$ is combinatorially free.

\section{A conjecture on the structure of $\pi_*(G(A))$}
\noindent{\bf Example 3: }Let $A$ be a rank $n$ Cartan matrix $A$ which satisfy $a_{ij}a_{ji}\geq 4$ for all $i\not =j$, then the Weyl group of $G(A)$ is
$$W(A)=<\sigma_1,\cdots,\sigma_n|\sigma^2_i=1,1\leq i\leq n>.$$

\noindent By Lemma 2, $\displaystyle{\frac{1}{H_A}=\frac{n}{1+q}-(n-1)}$, so $H_A=\displaystyle{\frac{1+q}{1-(n-1)q}}$.
$$C_A(q)=\frac{(1-q)^{n-1}}{1-(n-1)q}=\displaystyle{\frac{1}{\frac{1-(n-1)q}{(1-q)^{n-1}}}}= \frac{1}{1-a_2 q^2-a_3 q^3-\cdots-a_k q^k-\cdots}.$$ where $\displaystyle{a_k=(n-1) {{k+n-3}\choose {k-1}}-{{k+n-2}\choose {k}}}$.

\noindent{\bf Lemma 11: }$a_k> 0$ for all $k>1$.

\noindent{\bf Conjecture : }$\pi_{even}(G(A))$ is a free Lie algebra with $a_k$ free generators of degree $2k$ for each $k>1$.

%So by compare coefficients we have

%\noindent {\bf Lemma : }$c_2=k$. $c_4=?, c_6$.

%THe Lemma to assert a strongly free set in $L$

%On the degree problem for $H_A$ and $P_A$.

\newpage
\section{Appendix-the lists of Poincar\'{e} series of rank $3$ and $4$}

\begin{center}
{\bf The Poincar\'{e} series of Cartan matrices with reduced Coxeter graph 3-I}

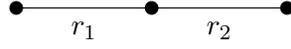
\begin{center}

\end{center}
\begin{tabular}{|c|c|c|}
\hline
& $[r1,r2]$ & $Poincar\acute{e}\, series$ \\\hline
1 & $[1,1]$ & $\displaystyle{\frac { \left( {t}^{2}-1 \right)  \left( {t}^{3}-1
 \right)  \left( {t}^{4}-1 \right) }{ \left( t-1 \right) ^{3}}}$\\\hline
2 & $[1,2]$ & $\displaystyle{\frac { \left( {t}^{2}-1 \right)  \left( {t}^{4}-1
 \right)  \left( {t}^{6}-1 \right) }{ \left( t-1 \right) ^{3}}}$\\\hline
3 & $[1,3]$ & $\displaystyle{\frac { \left( {t}^{5}+{t}^{4}+{t}^{3}+{t}^{2}+t+1
 \right)  \left( t+1 \right) }{{t}^{6}-{t}^{5}-t+1}}$\\\hline
4 & $[1,4]$ & $\displaystyle{-\frac { \left( {t}^{2}+t+1 \right)  \left( t+1 \right) ^{
2}}{{t}^{3}+{t}^{2}-1}}$\\\hline
5 & $[2,2]$ & $\displaystyle{\frac { \left( {t}^{3}+{t}^{2}+t+1 \right)  \left( t+1
 \right) }{{t}^{4}-{t}^{3}-t+1}}$\\\hline
6 & $[2,3]$ & $\displaystyle{\frac { \left( {t}^{5}+{t}^{4}+{t}^{3}+{t}^{2}+t+1
 \right)  \left( {t}^{3}+{t}^{2}+t+1 \right) }{{t}^{8}-{t}^{7}+{t}^{6}
-2\,{t}^{5}+{t}^{4}-2\,{t}^{3}+{t}^{2}-t+1}}$\\\hline
7 & $[2,4]$ & $\displaystyle{-\frac { \left( {t}^{3}+{t}^{2}+t+1 \right)  \left( t+1
 \right) }{{t}^{3}+t-1}}$\\\hline
8 & $[3,3]$ & $\displaystyle{\frac { \left( {t}^{5}+{t}^{4}+{t}^{3}+{t}^{2}+t+1
 \right)  \left( t+1 \right) }{{t}^{6}-{t}^{5}-{t}^{3}-t+1}}$\\\hline
9 & $[3,4]$ & $\displaystyle{-\frac { \left( {t}^{5}+{t}^{4}+{t}^{3}+{t}^{2}+t+1
 \right)  \left( t+1 \right) }{{t}^{5}+{t}^{3}+t-1}}$\\\hline
10 & $[4,4]$ & $\displaystyle{-\frac { \left( t+1 \right) ^{2}}{{t}^{2}+t-1}}$\\\hline
\end{tabular}
\end{center}

\newpage

\begin{center}
{\bf The Poincar\'{e} series of Cartan matrices with reduced Coxeter graph 3-II}



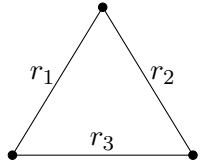
\begin{tabular}{|c|c|c|}
\hline
& $[r1,r2,r3]$ & $Poincar\acute{e}\, series$ \\\hline
11& $[1,1,1]$ & $\displaystyle{\frac {{t}^{2}+t+1}{{t}^{2}-2\,t+1}}$\\\hline
12& $[1,1,2]$ & $\displaystyle{\frac { \left( {t}^{3}+{t}^{2}+t+1 \right)  \left( {t}^{
2}+t+1 \right)  \left( t+1 \right) }{{t}^{6}-{t}^{4}-{t}^{3}-{t}^{2}+1
}}$\\\hline
13& $[1,1,3]$ & $\displaystyle{\frac { \left( {t}^{5}+{t}^{4}+{t}^{3}+{t}^{2}+t+1
 \right)  \left( t+1 \right) }{{t}^{6}-{t}^{5}-{t}^{4}+{t}^{3}-{t}^{2}
-t+1}}$\\\hline
14& $[1,1,4]$ & $\displaystyle{-\frac { \left( {t}^{2}+t+1 \right)  \left( t+1 \right)
}{{t}^{2}+t-1}}$\\\hline
15& $[1,2,2]$ & $\displaystyle{\frac { \left( {t}^{3}+{t}^{2}+t+1 \right)  \left( {t}^{
2}+t+1 \right)  \left( t+1 \right) }{{t}^{6}-{t}^{4}-2\,{t}^{3}-{t}^{2
}+1}}$\\\hline
16& $[1,2,3]$ & $\displaystyle{\frac { \left( {t}^{5}+{t}^{4}+{t}^{3}+{t}^{2}+t+1
 \right)  \left( {t}^{3}+{t}^{2}+t+1 \right) }{{t}^{8}-{t}^{7}-{t}^{5}
-{t}^{4}-{t}^{3}-t+1}}$\\\hline
17& $[1,2,4]$ & $\displaystyle{-\frac { \left( {t}^{3}+{t}^{2}+t+1 \right)  \left( {t}^
{2}+t+1 \right)  \left( t+1 \right) }{{t}^{5}+2\,{t}^{4}+2\,{t}^{3}+{t
}^{2}-1}}$\\
\hline
\end{tabular}
\end{center}
\newpage

\begin{center}
{\bf The Poincar\'{e} series of Cartan matrices with reduced Coxeter graph 4-I}
\end{center}

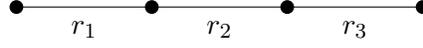
\begin{center}

\end{center}

\begin{center}
\begin{tabular}{|c|c|c|}
\hline
 & $[r1,r2,r3]$ & $Poincar\acute{e}\,\, series$ \\\hline
1 & $[1,1,1]$ & $\displaystyle{\frac { \left( {t}^{2}-1 \right)  \left( {t}^{3}-1
 \right)  \left( {t}^{4}-1 \right)  \left( {t}^{5}-1 \right) }{
 \left( t-1 \right) ^{4}}}$\\\hline
2 & $[1,1,2]$ & $\displaystyle{\frac { \left( {t}^{2}-1 \right)  \left( {t}^{4}-1
 \right)  \left( {t}^{6}-1 \right)  \left( {t}^{8}-1 \right) }{
 \left( t-1 \right) ^{4}}}$\\\hline
3 & $[1,1,3]$ & $\displaystyle{\frac { \left( {t}^{5}+{t}^{4}+{t}^{3}+{t}^{2}+t+1
 \right)  \left( {t}^{3}+{t}^{2}+t+1 \right)  \left( t+1 \right) }{{t}
^{8}-{t}^{4}-t+1}}$\\\hline
4 & $[1,1,4]$ & $\displaystyle{-\frac { \left( {t}^{2}+t+1 \right)  \left( {t}^{3}+{t}^
{2}+t+1 \right)  \left( t+1 \right) }{{t}^{3}+t-1}}$\\\hline
5& $[1,2,1]$ & $\displaystyle{\frac { \left( {t}^{2}-1 \right)  \left( {t}^{6}-1
 \right)  \left( {t}^{8}-1 \right)  \left( {t}^{12}-1 \right) }{
 \left( t-1 \right) ^{4}}}$\\\hline
6& $[1,2,2]$ & $\displaystyle{\frac { \left( {t}^{5}+{t}^{4}+{t}^{3}+{t}^{2}+t+1
 \right)  \left( {t}^{3}+{t}^{2}+t+1 \right)  \left( t+1 \right) }{{t}
^{8}-{t}^{3}-t+1}}$\\\hline
7& $[1,2,3]$ & $\displaystyle{\frac { \left( {t}^{5}+{t}^{4}+{t}^{3}+{t}^{2}+t+1
 \right)  \left( {t}^{3}+{t}^{2}+t+1 \right)  \left( t+1 \right) }{{t}
^{8}-{t}^{4}-{t}^{3}-t+1}}$\\
\hline

\end{tabular}

\begin{tabular}{|c|c|c|}
\hline
 & $[r1,r2,r3]$ & $Poincar\acute{e}\,\, series$ \\\hline

21& $[2,2,3]$ & $\displaystyle{\frac { \left( {t}^{3}+{t}^{2}+t+1 \right)  \left( {t}^
{5}+{t}^{4}+{t}^{3}+{t}^{2}+t+1 \right)  \left( t+1 \right) }{{t}^{9}+
{t}^{8}-{t}^{5}-2\,{t}^{3}-t+1}}$\\\hline
22& $[2,2,4]$ & $\displaystyle{\frac { \left( {t}^{3}+{t}^{2}+t+1 \right)  \left( t+1
 \right) ^{2}}{{t}^{5}-{t}^{3}-{t}^{2}-t+1}}$\\\hline
23& $[2,3,2]$ & $\displaystyle{\frac { \left( {t}^{3}+{t}^{2}+t+1 \right)  \left( {t}^
{5}+{t}^{4}+{t}^{3}+{t}^{2}+t+1 \right)  \left( t+1 \right) }{{t}^{9}+
{t}^{8}-2\,{t}^{5}-2\,{t}^{3}-t+1}}$\\\hline
24& $[2,3,3]$ & $\displaystyle{\frac { \left( {t}^{3}+{t}^{2}+t+1 \right)  \left( {t}^
{5}+{t}^{4}+{t}^{3}+{t}^{2}+t+1 \right)  \left( t+1 \right) }{{t}^{9}+
{t}^{8}-2\,{t}^{5}-{t}^{4}-2\,{t}^{3}-t+1}}$\\\hline
25& $[2,3,4]$ & $\displaystyle{\frac { \left( {t}^{3}+{t}^{2}+t+1 \right)  \left( {t}^
{5}+{t}^{4}+{t}^{3}+{t}^{2}+t+1 \right)  \left( t+1 \right) }{{t}^{9}-
{t}^{6}-2\,{t}^{5}-{t}^{4}-2\,{t}^{3}-t+1}}$\\\hline
26& $[2,4,2]$ & $\displaystyle{-\frac { \left( {t}^{3}+{t}^{2}+t+1 \right)  \left( t+1
 \right) ^{2}}{{t}^{3}+{t}^{2}+t-1}}$\\\hline
27& $[2,4,3]$ & $\displaystyle{-\frac { \left( {t}^{5}+{t}^{4}+{t}^{3}+{t}^{2}+t+1
 \right)  \left( {t}^{3}+{t}^{2}+t+1 \right)  \left( t+1 \right) }{{t}
^{7}+{t}^{6}+2\,{t}^{5}+{t}^{4}+2\,{t}^{3}+t-1}}$\\
\hline
\end{tabular}

\newpage
\begin{center}
{\bf The Poincar\'{e} series of Cartan matrices with reduced Coxeter graph 4-II}
\end{center}

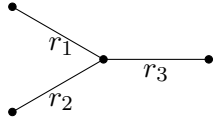
\begin{center}

\end{center}

\begin{tabular}{|c|c|c|}
\hline
 & $[r1,r2,r3]$ & $Poincar\acute{e}\,\, series$ \\\hline
41 & $[1,1,1]$ & $\displaystyle{\frac { \left( {t}^{2}-1 \right)  \left( {t}^{4}-1
 \right) ^{2} \left( {t}^{6}-1 \right) }{ \left( t-1 \right) ^{4}}}$\\\hline
42 & $[1,1,2]$ & $\displaystyle{-\frac { \left( {t}^{3}+1 \right)  \left( {t}^{3}+{t}^{
2}+t+1 \right)  \left( t+1 \right) }{{t}^{7}-2\,{t}^{6}+{t}^{5}-{t}^{2
}+2\,t-1}}$\\\hline
43 & $[1,1,3]$ & $\displaystyle{\frac { \left( {t}^{5}+{t}^{4}+{t}^{3}+{t}^{2}+t+1
 \right)  \left( {t}^{3}+{t}^{2}+t+1 \right)  \left( t+1 \right) }{{t}
^{9}+{t}^{8}+{t}^{6}-2\,{t}^{5}-{t}^{3}-t+1}}$\\\hline
44 & $[1,1,4]$ & $\displaystyle{-\frac { \left( {t}^{2}+t+1 \right)  \left( {t}^{3}+{t}
^{2}+t+1 \right)  \left( t+1 \right) ^{2}}{{t}^{5}+{t}^{4}+2\,{t}^{3}+
{t}^{2}-1}}$\\\hline
45 & $[1,2,2]$ & $\displaystyle{\frac { \left( {t}^{3}+1 \right)  \left( {t}^{3}+{t}^{2
}+t+1 \right)  \left( t+1 \right) }{{t}^{6}-{t}^{5}+{t}^{4}-{t}^{3}+{t
}^{2}-2\,t+1}}$\\\hline
46 & $[1,2,3]$ & $\displaystyle{\frac { \left( {t}^{5}+{t}^{4}+{t}^{3}+{t}^{2}+t+1
 \right)  \left( {t}^{3}+{t}^{2}+t+1 \right)  \left( t+1 \right) }{{t}
^{9}+{t}^{8}+{t}^{6}-2\,{t}^{5}-2\,{t}^{3}-t+1}}$\\\hline
47 & $[1,2,4]$ & $\displaystyle{-\frac { \left( {t}^{2}+t+1 \right)  \left( {t}^{3}+1
 \right)  \left( {t}^{3}+{t}^{2}+t+1 \right)  \left( t+1 \right) }{{t}
^{7}+2\,{t}^{5}+2\,{t}^{3}+t-1}}$\\
\hline
\end{tabular}

\begin{center}
{\bf The Poincar\'{e} series of Cartan matrices with reduced Coxeter graph 4-III}
\end{center}

\begin{center}

\end{center}

\begin{tabular}{|c|c|c|c|}
\hline
 & $[r1,r2,r3,r4]$ & $Poincar\acute{e}\,\, series$ \\\hline

61 & $[1,1,1,1]$ & $\displaystyle{\frac { \left( {t}^{2}+t+1 \right)  \left( {t}^{3}+{t
}^{2}+t+1 \right) }{{t}^{4}-{t}^{3}+{t}^{2}-2\,t+1}}$\\\hline
62 & $[1,1,1,2]$ & $\displaystyle{\frac { \left( {t}^{2}+t+1 \right)  \left( {t}^{3}+{t
}^{2}+t+1 \right)  \left( t+1 \right) ^{2}}{{t}^{6}+{t}^{5}-2\,{t}^{3}
-2\,{t}^{2}+1}}$\\\hline
63 & $[1,1,1,3]$ & $\displaystyle{\frac { \left( {t}^{5}+{t}^{4}+{t}^{3}+{t}^{2}+t+1
 \right)  \left( {t}^{3}+{t}^{2}+t+1 \right)  \left( t+1 \right) }{{t}
^{8}-{t}^{5}-{t}^{4}-{t}^{2}-t+1}}$\\\hline
64 & $[1,1,1,4]$ & $\displaystyle{-\frac { \left( {t}^{3}+{t}^{2}+t+1 \right)  \left( {
t}^{2}+t+1 \right)  \left( t+1 \right) }{{t}^{3}+{t}^{2}+t-1}}$\\\hline
65 & $[1,1,2,1]$ & $\displaystyle{\frac { \left( {t}^{5}+{t}^{4}+{t}^{3}+{t}^{2}+t+1
 \right)  \left( {t}^{3}+{t}^{2}+t+1 \right)  \left( t+1 \right) }{{t}
^{8}-{t}^{4}-{t}^{2}-t+1}}$\\\hline
66 & $[1,1,2,2]$ & $\displaystyle{\frac { \left( {t}^{5}+{t}^{4}+{t}^{3}+{t}^{2}+t+1
 \right)  \left( {t}^{3}+{t}^{2}+t+1 \right)  \left( t+1 \right) }{{t}
^{8}-{t}^{5}-{t}^{3}-{t}^{2}-t+1}}$\\\hline
67 & $[1,1,2,3]$ & $\displaystyle{\frac { \left( {t}^{5}+{t}^{4}+{t}^{3}+{t}^{2}+t+1
 \right)  \left( {t}^{3}+{t}^{2}+t+1 \right)  \left( t+1 \right) }{{t}
^{8}-{t}^{5}-{t}^{4}-{t}^{3}-{t}^{2}-t+1}}$\\
\hline

\end{tabular}

\begin{tabular}{|c|c|c|c|}
\hline
 & $[r1,r2,r3,r4]$ & $Poincar\acute{e}\,\, series$ \\\hline

81 & $[1,2,3,1]$ & $\displaystyle{\frac { \left( {t}^{5}+{t}^{4}+{t}^{3}+{t}^{2}+t+1
 \right)  \left( {t}^{3}+{t}^{2}+t+1 \right)  \left( t+1 \right) }{{t}
^{9}+{t}^{8}-{t}^{5}-2\,{t}^{4}-{t}^{3}-{t}^{2}-t+1}}$\\\hline
82 & $[1,2,3,2]$ & $\displaystyle{\frac { \left( {t}^{5}+{t}^{4}+{t}^{3}+{t}^{2}+t+1
 \right)  \left( {t}^{3}+{t}^{2}+t+1 \right)  \left( t+1 \right) }{{t}
^{9}+{t}^{8}-2\,{t}^{5}-{t}^{4}-2\,{t}^{3}-{t}^{2}-t+1}}$\\\hline
83 & $[1,2,3,3]$ & $\displaystyle{\frac { \left( {t}^{5}+{t}^{4}+{t}^{3}+{t}^{2}+t+1
 \right)  \left( {t}^{3}+{t}^{2}+t+1 \right)  \left( t+1 \right) }{{t}
^{9}+{t}^{8}-2\,{t}^{5}-2\,{t}^{4}-2\,{t}^{3}-{t}^{2}-t+1}}$\\\hline
84 & $[1,2,3,4]$ & $\displaystyle{\frac { \left( {t}^{5}+{t}^{4}+{t}^{3}+{t}^{2}+t+1
 \right)  \left( {t}^{3}+{t}^{2}+t+1 \right)  \left( t+1 \right) }{{t}
^{9}-{t}^{6}-2\,{t}^{5}-2\,{t}^{4}-2\,{t}^{3}-{t}^{2}-t+1}}$\\\hline
85 & $[1,2,4,1]$ & $\displaystyle{-\frac { \left( {t}^{2}+t+1 \right)  \left( {t}^{3}+1
 \right)  \left( {t}^{3}+{t}^{2}+t+1 \right)  \left( t+1 \right) }{{t}
^{7}+{t}^{6}+{t}^{5}+2\,{t}^{4}+{t}^{3}+{t}^{2}+t-1}}$\\\hline
86 & $[1,2,4,2]$ & $\displaystyle{-\frac { \left( {t}^{2}+t+1 \right)  \left( {t}^{3}+{
t}^{2}+t+1 \right)  \left( {t}^{3}+1 \right)  \left( t+1 \right) }{{t}
^{7}+{t}^{6}+2\,{t}^{5}+{t}^{4}+2\,{t}^{3}+{t}^{2}+t-1}}$\\\hline
87 & $[1,2,4,3]$ & $\displaystyle{-\frac { \left( {t}^{2}+t+1 \right)  \left( {t}^{3}+{
t}^{2}+t+1 \right)  \left( {t}^{3}+1 \right)  \left( t+1 \right) }{{t}
^{7}+{t}^{6}+2\,{t}^{5}+2\,{t}^{4}+2\,{t}^{3}+{t}^{2}+t-1}}$\\
\hline

\end{tabular}

\newpage

\begin{tabular}{|c|c|c|c|}
\hline
 & $[r1,r2,r3,r4]$ & $Poincar\acute{e}\,\, series$ \\\hline

101 & $[2,1,1,1]$ & $\displaystyle{\frac { \left( {t}^{2}+t+1 \right)  \left( {t}^{3}+1
 \right)  \left( {t}^{3}+{t}^{2}+t+1 \right)  \left( t+1 \right) }{{t}
^{8}+{t}^{6}-{t}^{4}-{t}^{2}-t+1}}$\\\hline
102 & $[2,1,1,2]$ & $\displaystyle{\frac { \left( {t}^{2}+t+1 \right)  \left( {t}^{3}+1
 \right)  \left( {t}^{3}+{t}^{2}+t+1 \right)  \left( t+1 \right) }{{t}
^{8}+{t}^{6}-{t}^{5}-{t}^{3}-{t}^{2}-t+1}}$\\\hline
103 & $[2,1,1,3]$ & $\displaystyle{\frac { \left( {t}^{5}+{t}^{4}+{t}^{3}+{t}^{2}+t+1
 \right)  \left( {t}^{3}+{t}^{2}+t+1 \right)  \left( t+1 \right) }{{t}
^{8}+{t}^{6}-{t}^{5}-{t}^{4}-{t}^{3}-{t}^{2}-t+1}}$\\\hline
104 & $[2,1,1,4]$ & $\displaystyle{-\frac { \left( {t}^{2}+t+1 \right)  \left( {t}^{3}+
1 \right)  \left( {t}^{3}+{t}^{2}+t+1 \right)  \left( t+1 \right) }{{t
}^{5}+{t}^{4}+{t}^{3}+{t}^{2}+t-1}}$\\\hline
105 & $[2,1,2,1]$ & $\displaystyle{\frac { \left( {t}^{2}+t+1 \right)  \left( {t}^{3}+1
 \right)  \left( {t}^{3}+{t}^{2}+t+1 \right)  \left( t+1 \right) }{{t}
^{9}+{t}^{8}-{t}^{4}-{t}^{3}-{t}^{2}-t+1}}$\\\hline
106 & $[2,1,2,2]$ & $\displaystyle{\frac { \left( {t}^{2}+t+1 \right)  \left( {t}^{3}+1
 \right)  \left( {t}^{3}+{t}^{2}+t+1 \right)  \left( t+1 \right) }{{t}
^{9}+{t}^{8}-{t}^{5}-2\,{t}^{3}-{t}^{2}-t+1}}$\\\hline
107 & $[2,1,2,3]$ & $\displaystyle{\frac { \left( {t}^{5}+{t}^{4}+{t}^{3}+{t}^{2}+t+1
 \right)  \left( {t}^{3}+{t}^{2}+t+1 \right)  \left( t+1 \right) }{{t}
^{9}+{t}^{8}-{t}^{5}-{t}^{4}-2\,{t}^{3}-{t}^{2}-t+1}}$\\
\hline

\end{tabular}

\begin{tabular}{|c|c|c|c|}
\hline
 & $[r1,r2,r3,r4]$ & $Poincar\acute{e}\,\, series$ \\\hline

121 & $[2,2,3,1]$ & $\displaystyle{\frac { \left( {t}^{3}+{t}^{2}+t+1 \right)  \left( {
t}^{5}+{t}^{4}+{t}^{3}+{t}^{2}+t+1 \right) }{2\,{t}^{8}-{t}^{7}+{t}^{6
}-2\,{t}^{5}+{t}^{4}-3\,{t}^{3}+{t}^{2}-2\,t+1}}$\\\hline
122 & $[2,2,3,2]$ & $\displaystyle{\frac { \left( {t}^{3}+{t}^{2}+t+1 \right)  \left( {
t}^{5}+{t}^{4}+{t}^{3}+{t}^{2}+t+1 \right)  \left( t+1 \right) }{2\,{t
}^{9}+{t}^{8}-{t}^{6}-2\,{t}^{5}-{t}^{4}-3\,{t}^{3}-{t}^{2}-t+1}}$\\\hline
123 & $[2,2,3,3]$ & $\displaystyle{\frac { \left( {t}^{3}+{t}^{2}+t+1 \right)  \left( {
t}^{5}+{t}^{4}+{t}^{3}+{t}^{2}+t+1 \right)  \left( t+1 \right) }{2\,{t
}^{9}+{t}^{8}-{t}^{6}-2\,{t}^{5}-2\,{t}^{4}-3\,{t}^{3}-{t}^{2}-t+1}}$\\\hline
124 & $[2,2,3,4]$ & $\displaystyle{\frac { \left( {t}^{3}+{t}^{2}+t+1 \right)  \left( {
t}^{5}+{t}^{4}+{t}^{3}+{t}^{2}+t+1 \right) }{2\,{t}^{8}-2\,{t}^{7}+2\,
{t}^{6}-4\,{t}^{5}+2\,{t}^{4}-4\,{t}^{3}+{t}^{2}-2\,t+1}}$\\\hline
125 & $[2,2,4,1]$ & $\displaystyle{\frac { \left( {t}^{3}+{t}^{2}+t+1 \right)  \left( {
t}^{2}+t+1 \right)  \left( t+1 \right) }{{t}^{6}-{t}^{4}-3\,{t}^{3}-{t
}^{2}-t+1}}$\\\hline
126 & $[2,2,4,2]$ & $\displaystyle{\frac { \left( {t}^{3}+{t}^{2}+t+1 \right)  \left( t
+1 \right) ^{2}}{{t}^{5}-2\,{t}^{3}-2\,{t}^{2}-t+1}}$\\\hline
127 & $[2,2,4,3]$ & $\displaystyle{\frac { \left( {t}^{3}+{t}^{2}+t+1 \right)  \left( {
t}^{5}+{t}^{4}+{t}^{3}+{t}^{2}+t+1 \right)  \left( t+1 \right) }{{t}^{
9}-{t}^{7}-2\,{t}^{6}-2\,{t}^{5}-2\,{t}^{4}-3\,{t}^{3}-{t}^{2}-t+1}}$\\
\hline

\end{tabular}

\begin{tabular}{|c|c|c|c|}
\hline
 & $[r1,r2,r3,r4]$ & $Poincar\acute{e}\,\, series$ \\\hline

141 & $[3,1,1,1]$ & $\displaystyle{\frac { \left( {t}^{2}+t+1 \right)  \left( {t}^{3}+1
 \right)  \left( t+1 \right) }{2\,{t}^{6}-{t}^{5}-{t}^{4}+{t}^{3}-2\,t
+1}}$\\\hline
142 & $[3,1,1,2]$ & $\displaystyle{\frac { \left( {t}^{3}+{t}^{2}+t+1 \right)  \left( {
t}^{2}+t+1 \right)  \left( {t}^{3}+1 \right)  \left( t+1 \right) }{2\,
{t}^{9}+{t}^{8}+{t}^{6}-2\,{t}^{5}-{t}^{4}-{t}^{3}-{t}^{2}-t+1}}$\\\hline
143 & $[3,1,1,3]$ & $\displaystyle{\frac { \left( {t}^{3}+1 \right)  \left( {t}^{2}+t+1
 \right)  \left( t+1 \right) ^{2}}{2\,{t}^{7}+{t}^{6}-2\,{t}^{5}-2\,{t
}^{2}-t+1}}$\\\hline
144 & $[3,1,1,4]$ & $\displaystyle{\frac { \left( {t}^{2}+t+1 \right)  \left( {t}^{3}+1
 \right)  \left( t+1 \right) }{2\,{t}^{6}-2\,{t}^{5}-2\,t+1}}$\\\hline
145 & $[3,1,2,1]$ & $\displaystyle{\frac { \left( {t}^{5}+{t}^{4}+{t}^{3}+{t}^{2}+t+1
 \right)  \left( {t}^{3}+{t}^{2}+t+1 \right) }{2\,{t}^{8}-{t}^{7}+{t}^
{6}-{t}^{5}-2\,{t}^{3}+{t}^{2}-2\,t+1}}$\\\hline
146 & $[3,1,2,2]$ & $\displaystyle{\frac { \left( {t}^{5}+{t}^{4}+{t}^{3}+{t}^{2}+t+1
 \right)  \left( {t}^{3}+{t}^{2}+t+1 \right)  \left( t+1 \right) }{2\,
{t}^{9}+{t}^{8}-2\,{t}^{5}-{t}^{4}-2\,{t}^{3}-{t}^{2}-t+1}}$\\\hline
147 & $[3,1,2,3]$ & $\displaystyle{\frac { \left( {t}^{5}+{t}^{4}+{t}^{3}+{t}^{2}+t+1
 \right)  \left( {t}^{3}+{t}^{2}+t+1 \right)  \left( t+1 \right) }{2\,
{t}^{9}+{t}^{8}-2\,{t}^{5}-2\,{t}^{4}-2\,{t}^{3}-{t}^{2}-t+1}}$\\
\hline

\end{tabular}

\begin{tabular}{|c|c|c|c|}
\hline
 & $[r1,r2,r3,r4]$ & $Poincar\acute{e}\,\, series$ \\\hline

161 & $[3,2,3,1]$ & $\displaystyle{\frac { \left( {t}^{5}+{t}^{4}+{t}^{3}+{t}^{2}+t+1
 \right)  \left( {t}^{3}+{t}^{2}+t+1 \right) }{2\,{t}^{8}-{t}^{7}+{t}^
{6}-2\,{t}^{5}-3\,{t}^{3}+{t}^{2}-2\,t+1}}$\\\hline
162 & $[3,2,3,2]$ & $\displaystyle{\frac { \left( {t}^{3}+{t}^{2}+t+1 \right)  \left( {
t}^{5}+{t}^{4}+{t}^{3}+{t}^{2}+t+1 \right)  \left( t+1 \right) }{2\,{t
}^{9}+{t}^{8}-{t}^{6}-3\,{t}^{5}-2\,{t}^{4}-3\,{t}^{3}-{t}^{2}-t+1}}$\\\hline
163 & $[3,2,3,3]$ & $\displaystyle{\frac { \left( {t}^{3}+{t}^{2}+t+1 \right)  \left( {
t}^{5}+{t}^{4}+{t}^{3}+{t}^{2}+t+1 \right)  \left( t+1 \right) }{2\,{t
}^{9}+{t}^{8}-{t}^{6}-3\,{t}^{5}-3\,{t}^{4}-3\,{t}^{3}-{t}^{2}-t+1}}$\\\hline
164 & $[3,2,3,4]$ & $\displaystyle{\frac { \left( {t}^{5}+{t}^{4}+{t}^{3}+{t}^{2}+t+1
 \right)  \left( {t}^{3}+{t}^{2}+t+1 \right) }{2\,{t}^{8}-2\,{t}^{7}+2
\,{t}^{6}-4\,{t}^{5}+{t}^{4}-4\,{t}^{3}+{t}^{2}-2\,t+1}}$\\\hline
165 & $[3,2,4,1]$ & $\displaystyle{\frac { \left( {t}^{2}+t+1 \right)  \left( {t}^{3}+{
t}^{2}+t+1 \right)  \left( {t}^{3}+1 \right) }{{t}^{8}-{t}^{7}-2\,{t}^
{5}-3\,{t}^{3}+{t}^{2}-2\,t+1}}$\\\hline
166 & $[3,2,4,2]$ & $\displaystyle{\frac { \left( {t}^{3}+{t}^{2}+t+1 \right)  \left( {
t}^{5}+{t}^{4}+{t}^{3}+{t}^{2}+t+1 \right)  \left( t+1 \right) }{{t}^{
9}-{t}^{7}-2\,{t}^{6}-3\,{t}^{5}-2\,{t}^{4}-3\,{t}^{3}-{t}^{2}-t+1}}$\\\hline
167 & $[3,2,4,3]$ & $\displaystyle{\frac { \left( {t}^{3}+{t}^{2}+t+1 \right)  \left( {
t}^{5}+{t}^{4}+{t}^{3}+{t}^{2}+t+1 \right)  \left( t+1 \right) }{{t}^{
9}-{t}^{7}-2\,{t}^{6}-3\,{t}^{5}-3\,{t}^{4}-3\,{t}^{3}-{t}^{2}-t+1}}$\\
\hline
\end{tabular}

\begin{tabular}{|c|c|c|c|}
\hline
 & $[r1,r2,r3,r4]$ & $Poincar\acute{e}\,\, series$ \\\hline

181 & $[4,1,1,1]$ & $\displaystyle{\frac { \left( {t}^{2}+t+1 \right)  \left( t+1
 \right) ^{2}}{{t}^{4}-2\,{t}^{2}-t+1}}$\\\hline
182 & $[4,1,1,2]$ & $\displaystyle{\frac { \left( {t}^{3}+{t}^{2}+t+1 \right)  \left( {
t}^{2}+t+1 \right)  \left( t+1 \right) ^{2}}{{t}^{7}+{t}^{6}-{t}^{5}-2
\,{t}^{4}-3\,{t}^{3}-2\,{t}^{2}+1}}$\\\hline
183 & $[4,1,1,3]$ & $\displaystyle{\frac { \left( {t}^{5}+{t}^{4}+{t}^{3}+{t}^{2}+t+1
 \right)  \left( t+1 \right) ^{2}}{{t}^{7}-2\,{t}^{5}-2\,{t}^{2}-t+1}}
]$\\\hline
184 & $[4,1,1,4]$ & $\displaystyle{\frac { \left( {t}^{2}+t+1 \right)  \left( t+1
 \right) ^{2}}{{t}^{4}-{t}^{3}-2\,{t}^{2}-t+1}}$\\\hline
185 & $[4,1,2,1]$ & $\displaystyle{\frac { \left( {t}^{3}+{t}^{2}+t+1 \right)  \left( {
t}^{2}+t+1 \right)  \left( t+1 \right) }{{t}^{6}-{t}^{4}-2\,{t}^{3}-{t
}^{2}-t+1}}$\\\hline
186 & $[4,1,2,2]$ & $\displaystyle{\frac { \left( {t}^{3}+{t}^{2}+t+1 \right)  \left( {
t}^{2}+t+1 \right)  \left( t+1 \right) ^{2}}{{t}^{7}+{t}^{6}-{t}^{5}-3
\,{t}^{4}-4\,{t}^{3}-2\,{t}^{2}+1}}$\\\hline
187 & $[4,1,2,3]$ & $\displaystyle{\frac { \left( {t}^{3}+{t}^{2}+t+1 \right)  \left( {
t}^{5}+{t}^{4}+{t}^{3}+{t}^{2}+t+1 \right)  \left( t+1 \right) }{{t}^{
9}-{t}^{7}-{t}^{6}-2\,{t}^{5}-2\,{t}^{4}-2\,{t}^{3}-{t}^{2}-t+1}}$\\
\hline

\end{tabular}

\begin{tabular}{|c|c|c|c|}
\hline
 & $[r1,r2,r3,r4]$ & $Poincar\acute{e}\,\, series$ \\\hline

201 & $[4,2,3,1]$ & $\displaystyle{\frac { \left( {t}^{5}+{t}^{4}+{t}^{3}+{t}^{2}+t+1
 \right)  \left( {t}^{3}+{t}^{2}+t+1 \right) }{{t}^{8}-{t}^{7}-2\,{t}^
{5}-3\,{t}^{3}+{t}^{2}-2\,t+1}}$\\\hline
202 & $[4,2,3,2]$ & $\displaystyle{\frac { \left( {t}^{5}+{t}^{4}+{t}^{3}+{t}^{2}+t+1
 \right)  \left( {t}^{3}+{t}^{2}+t+1 \right)  \left( t+1 \right) }{{t}
^{9}-{t}^{7}-2\,{t}^{6}-3\,{t}^{5}-2\,{t}^{4}-3\,{t}^{3}-{t}^{2}-t+1}}
]$\\\hline
203 & $[4,2,3,3]$ & $\displaystyle{\frac { \left( {t}^{5}+{t}^{4}+{t}^{3}+{t}^{2}+t+1
 \right)  \left( {t}^{3}+{t}^{2}+t+1 \right)  \left( t+1 \right) }{{t}
^{9}-{t}^{7}-2\,{t}^{6}-3\,{t}^{5}-3\,{t}^{4}-3\,{t}^{3}-{t}^{2}-t+1}}
]$\\\hline
204 & $[4,2,3,4]$ & $\displaystyle{\frac { \left( {t}^{5}+{t}^{4}+{t}^{3}+{t}^{2}+t+1
 \right)  \left( {t}^{3}+{t}^{2}+t+1 \right) }{{t}^{8}-2\,{t}^{7}+{t}^
{6}-4\,{t}^{5}+{t}^{4}-4\,{t}^{3}+{t}^{2}-2\,t+1}}$\\\hline
205 & $[4,2,4,1]$ & $\displaystyle{-\frac { \left( {t}^{3}+{t}^{2}+t+1 \right)  \left(
{t}^{2}+t+1 \right)  \left( t+1 \right) }{{t}^{5}+2\,{t}^{4}+3\,{t}^{3
}+{t}^{2}+t-1}}$\\\hline
206 & $[4,2,4,2]$ & $\displaystyle{-\frac { \left( {t}^{3}+{t}^{2}+t+1 \right)  \left(
t+1 \right) ^{2}}{{t}^{4}+2\,{t}^{3}+2\,{t}^{2}+t-1}}$\\\hline
207 & $[4,2,4,3]$ & $\displaystyle{-\frac { \left( {t}^{5}+{t}^{4}+{t}^{3}+{t}^{2}+t+1
 \right)  \left( {t}^{3}+{t}^{2}+t+1 \right)  \left( t+1 \right) }{{t}
^{8}+2\,{t}^{7}+3\,{t}^{6}+3\,{t}^{5}+3\,{t}^{4}+3\,{t}^{3}+{t}^{2}+t-
1}}$\\
\hline

\end{tabular}

\begin{center}
{\bf The Poincar\'{e} series of Cartan matrices with reduced Coxeter graph 4-IV}
\end{center}

\newpage
\begin{center}

\end{center}

\begin{tabular}{|c|c|c|c|}
\hline
 & $[r1,r2,r3,r4]$ & $Poincar\acute{e}\,\, series$ \\\hline

221 & $[1,1,1,1]$ & $\displaystyle{-\frac {{t}^{3}+{t}^{2}+t+1}{{t}^{3}-3\,{t}^{2}+3\,t
-1}}$\\\hline
222 & $[1,1,1,2]$ & $\displaystyle{-\frac { \left( {t}^{3}+1 \right)  \left( {t}^{3}+{t
}^{2}+t+1 \right)  \left( t+1 \right) }{{t}^{7}-2\,{t}^{6}+2\,{t}^{4}-
2\,{t}^{3}+2\,t-1}}$\\\hline
223 & $[1,1,1,3]$ & $\displaystyle{\frac { \left( {t}^{3}+1 \right)  \left( {t}^{2}+t+1
 \right)  \left( {t}^{3}+{t}^{2}+t+1 \right) }{{t}^{8}+{t}^{6}-{t}^{5}
-{t}^{3}+{t}^{2}-2\,t+1}}$\\\hline
224 & $[1,1,1,4]$ & $\displaystyle{-\frac { \left( {t}^{2}+t+1 \right)  \left( {t}^{3}+
{t}^{2}+t+1 \right)  \left( t+1 \right) }{{t}^{3}+{t}^{2}+t-1}}$\\\hline
225 & $[1,1,2,2]$ & $\displaystyle{\frac { \left( {t}^{3}+{t}^{2}+t+1 \right)  \left( {
t}^{3}+1 \right)  \left( t+1 \right) }{{t}^{6}-{t}^{4}+{t}^{3}-2\,t+1}
}$\\\hline
226 & $[1,1,2,3]$ & $\displaystyle{\frac { \left( {t}^{3}+1 \right)  \left( {t}^{2}+t+1
 \right)  \left( {t}^{3}+{t}^{2}+t+1 \right)  \left( t+1 \right) }{{t}
^{9}+{t}^{8}+{t}^{7}-{t}^{5}-{t}^{4}-{t}^{3}-{t}^{2}-t+1}}$\\\hline
227 & $[1,1,2,4]$ & $\displaystyle{-\frac { \left( {t}^{3}+1 \right)  \left( {t}^{2}+t+
1 \right)  \left( {t}^{3}+{t}^{2}+t+1 \right)  \left( t+1 \right) }{{t
}^{6}+{t}^{5}+{t}^{4}+{t}^{3}+{t}^{2}+t-1}}$\\
\hline

\end{tabular}

\begin{tabular}{|c|c|c|c|}
\hline
 & $[r1,r2,r3,r4]$ & $Poincar\acute{e}\,\, series$ \\\hline

241 & $[1,2,4,3]$ & $\displaystyle{\frac { \left( {t}^{2}+t+1 \right)  \left( {t}^{3}+{
t}^{2}+t+1 \right)  \left( {t}^{3}+1 \right)  \left( t+1 \right) }{{t}
^{9}-{t}^{6}-2\,{t}^{5}-2\,{t}^{4}-2\,{t}^{3}-{t}^{2}-t+1}}$\\\hline
242 & $[1,2,4,4]$ & $\displaystyle{-\frac { \left( {t}^{2}+t+1 \right)  \left( {t}^{3}+
{t}^{2}+t+1 \right)  \left( {t}^{3}+1 \right)  \left( t+1 \right) }{{t
}^{8}+{t}^{7}+2\,{t}^{6}+2\,{t}^{5}+2\,{t}^{4}+2\,{t}^{3}+{t}^{2}+t-1}
}$\\\hline
243 & $[1,3,1,3]$ & $\displaystyle{\frac { \left( {t}^{2}+t+1 \right)  \left( {t}^{3}+1
 \right)  \left( t+1 \right) }{3\,{t}^{6}-2\,{t}^{5}-2\,t+1}}$\\\hline
244 & $[1,3,1,4]$ & $\displaystyle{\frac { \left( {t}^{2}+t+1 \right)  \left( {t}^{3}+1
 \right)  \left( t+1 \right) }{2\,{t}^{6}-2\,{t}^{5}-2\,t+1}}$\\\hline
245 & $[1,3,2,3]$ & $\displaystyle{\frac { \left( {t}^{3}+1 \right)  \left( {t}^{2}+t+1
 \right)  \left( t+1 \right)  \left( {t}^{2}+1 \right) }{3\,{t}^{8}-2
\,{t}^{7}+3\,{t}^{6}-3\,{t}^{5}+{t}^{4}-3\,{t}^{3}+{t}^{2}-2\,t+1}}$\\\hline
246 & $[1,3,2,4]$ & $\displaystyle{\frac { \left( {t}^{3}+1 \right)  \left( {t}^{2}+t+1
 \right)  \left( t+1 \right)  \left( {t}^{2}+1 \right) }{2\,{t}^{8}-2
\,{t}^{7}+2\,{t}^{6}-3\,{t}^{5}+{t}^{4}-3\,{t}^{3}+{t}^{2}-2\,t+1}}$\\\hline
247 & $[1,3,3,3]$ & $\displaystyle{\frac { \left( {t}^{5}+{t}^{4}+{t}^{3}+{t}^{2}+t+1
 \right)  \left( t+1 \right) }{3\,{t}^{6}-2\,{t}^{5}-{t}^{3}-2\,t+1}}$\\
\hline

\end{tabular}

\begin{tabular}{|c|c|c|c|}
\hline
 & $[r1,r2,r3,r4]$ & $Poincar\acute{e}\,\, series$ \\\hline

261 & $[2,3,2,3]$ & $\displaystyle{\frac { \left( {t}^{3}+{t}^{2}+t+1 \right)  \left( {
t}^{5}+{t}^{4}+{t}^{3}+{t}^{2}+t+1 \right) }{3\,{t}^{8}-2\,{t}^{7}+3\,
{t}^{6}-4\,{t}^{5}+2\,{t}^{4}-4\,{t}^{3}+{t}^{2}-2\,t+1}}$\\\hline
262 & $[2,3,2,4]$ & $\displaystyle{\frac { \left( {t}^{3}+{t}^{2}+t+1 \right)  \left( {
t}^{5}+{t}^{4}+{t}^{3}+{t}^{2}+t+1 \right) }{2\,{t}^{8}-2\,{t}^{7}+2\,
{t}^{6}-4\,{t}^{5}+2\,{t}^{4}-4\,{t}^{3}+{t}^{2}-2\,t+1}}$\\\hline
263 & $[2,3,3,3]$ & $\displaystyle{\frac { \left( {t}^{3}+{t}^{2}+t+1 \right)  \left( {
t}^{5}+{t}^{4}+{t}^{3}+{t}^{2}+t+1 \right) }{3\,{t}^{8}-2\,{t}^{7}+3\,
{t}^{6}-4\,{t}^{5}+{t}^{4}-4\,{t}^{3}+{t}^{2}-2\,t+1}}$\\\hline
264 & $[2,3,3,4]$ & $\displaystyle{\frac { \left( {t}^{3}+{t}^{2}+t+1 \right)  \left( {
t}^{5}+{t}^{4}+{t}^{3}+{t}^{2}+t+1 \right) }{2\,{t}^{8}-2\,{t}^{7}+2\,
{t}^{6}-4\,{t}^{5}+{t}^{4}-4\,{t}^{3}+{t}^{2}-2\,t+1}}$\\\hline
265 & $[2,3,4,3]$ & $\displaystyle{\frac { \left( {t}^{3}+{t}^{2}+t+1 \right)  \left( {
t}^{5}+{t}^{4}+{t}^{3}+{t}^{2}+t+1 \right) }{2\,{t}^{8}-2\,{t}^{7}+2\,
{t}^{6}-4\,{t}^{5}+{t}^{4}-4\,{t}^{3}+{t}^{2}-2\,t+1}}$\\\hline
266 & $[2,3,4,4]$ & $\displaystyle{\frac { \left( {t}^{3}+{t}^{2}+t+1 \right)  \left( {
t}^{5}+{t}^{4}+{t}^{3}+{t}^{2}+t+1 \right) }{{t}^{8}-2\,{t}^{7}+{t}^{6
}-4\,{t}^{5}+{t}^{4}-4\,{t}^{3}+{t}^{2}-2\,t+1}}$\\\hline
267 & $[2,4,2,4]$ & $\displaystyle{\frac { \left( {t}^{3}+{t}^{2}+t+1 \right)  \left( t
+1 \right) }{{t}^{4}-2\,{t}^{3}-2\,t+1}}$\\
\hline

\end{tabular}

\newpage

\begin{center}
{\bf The Poincar\'{e} series of Cartan matrices with reduced Coxeter graph 4-V}
\end{center}

\begin{center}

\end{center}

\begin{tabular}{|c|c|c|c|c|}
\hline
 & $[r1,r2,r3,r4,r5]$ & $Poincar\acute{e}\,\, series$ \\\hline

276 & $[1,1,1,1,1]$ & $\displaystyle{\frac { \left( {t}^{2}+t+1 \right)
\left( {t}^{3}
+{t}^{2}+t+1 \right) }{{t}^{5}-2\,t+1}}$\\\hline
277 & $[1,1,1,1,2]$ & $\displaystyle{\frac { \left( {t}^{3}+{t}^{2}+t+1
\right)
 \left( {t}^{2}+t+1 \right)  \left( t+1 \right) }{{t}^{6}+{t}^{5}-{t}^
{3}-2\,{t}^{2}-t+1}}$\\\hline
278 & $[1,1,1,1,3]$ & $\displaystyle{\frac { \left( {t}^{5}+{t}^{4}+{t}^
{3}+{t}^{2}+t+1
 \right)  \left( {t}^{3}+{t}^{2}+t+1 \right) }{{t}^{8}-2\,{t}^{4}-2\,t
+1}}$\\\hline
279 & $[1,1,1,1,4]$ & $\displaystyle{-\frac { \left( {t}^{2}+t+1 \right)
\left( {t}^{3
}+{t}^{2}+t+1 \right) }{{t}^{3}+2\,t-1}}$\\\hline
280 & $[1,1,1,2,1]$ & $\displaystyle{\frac { \left( {t}^{3}+{t}^{2}+t+1
\right)
 \left( {t}^{5}+{t}^{4}+{t}^{3}+{t}^{2}+t+1 \right)  \left( t+1
 \right) }{{t}^{9}+{t}^{8}+{t}^{6}-{t}^{5}-{t}^{4}-2\,{t}^{2}-t+1}}$\\\hline
281 & $[1,1,1,2,2]$ & $\displaystyle{\frac { \left( {t}^{3}+{t}^{2}+t+1
\right)
 \left( {t}^{5}+{t}^{4}+{t}^{3}+{t}^{2}+t+1 \right)  \left( t+1
 \right) }{{t}^{9}+{t}^{8}-{t}^{5}-{t}^{4}-{t}^{3}-2\,{t}^{2}-t+1}}$\\\hline
282 & $[1,1,1,2,3]$ & $\displaystyle{\frac { \left( {t}^{3}+{t}^{2}+t+1
\right)
 \left( {t}^{5}+{t}^{4}+{t}^{3}+{t}^{2}+t+1 \right)  \left( t+1
 \right) }{{t}^{9}+{t}^{8}-2\,{t}^{5}-2\,{t}^{4}-{t}^{3}-2\,{t}^{2}-t+
1}}$\\
\hline

\end{tabular}

\begin{tabular}{|c|c|c|c|c|}
\hline
 & $[r1,r2,r3,r4,r5]$ & $Poincar\acute{e}\,\, series$ \\\hline

291 & $[1,1,1,4,4]$ & $\displaystyle{-\frac { \left( {t}^{2}+t+1 \right)
\left( {t}^{3
}+{t}^{2}+t+1 \right)  \left( t+1 \right) }{{t}^{5}+2\,{t}^{4}+2\,{t}^
{3}+2\,{t}^{2}+t-1}}$\\\hline
292 & $[1,1,2,2,1]$ & $\displaystyle{\frac { \left( {t}^{2}+t+1 \right)
\left( {t}^{3}
+{t}^{2}+t+1 \right)  \left( t+1 \right) }{2\,{t}^{6}+{t}^{5}-2\,{t}^{
3}-2\,{t}^{2}-t+1}}$\\\hline
293 & $[1,1,2,2,2]$ & $\displaystyle{\frac { \left( {t}^{2}+t+1 \right)
\left( {t}^{3}
+{t}^{2}+t+1 \right)  \left( t+1 \right) }{2\,{t}^{6}+{t}^{5}-3\,{t}^{
3}-2\,{t}^{2}-t+1}}$\\\hline
294 & $[1,1,2,2,3]$ & $\displaystyle{\frac { \left( {t}^{5}+{t}^{4}+{t}^
{3}+{t}^{2}+t+1
 \right)  \left( {t}^{3}+{t}^{2}+t+1 \right) }{2\,{t}^{8}-{t}^{7}+{t}^
{6}-2\,{t}^{5}-2\,{t}^{3}-2\,t+1}}$\\\hline
295 & $[1,1,2,2,4]$ & $\displaystyle{\frac { \left( {t}^{2}+t+1 \right)
\left( {t}^{3}
+{t}^{2}+t+1 \right)  \left( t+1 \right) }{{t}^{6}-{t}^{4}-3\,{t}^{3}-
2\,{t}^{2}-t+1}}$\\\hline
296 & $[1,1,2,3,1]$ & $\displaystyle{\frac { \left( {t}^{5}+{t}^{4}+{t}^
{3}+{t}^{2}+t+1
 \right)  \left( {t}^{3}+{t}^{2}+t+1 \right) }{2\,{t}^{8}-{t}^{7}+{t}^
{6}-{t}^{5}-{t}^{4}-{t}^{3}-2\,t+1}}$\\\hline
297 & $[1,1,2,3,2]$ & $\displaystyle{\frac { \left( {t}^{5}+{t}^{4}+{t}^
{3}+{t}^{2}+t+1
 \right)  \left( {t}^{3}+{t}^{2}+t+1 \right) }{2\,{t}^{8}-{t}^{7}+{t}^
{6}-2\,{t}^{5}-2\,{t}^{3}-2\,t+1}}$\\
\hline

\end{tabular}

\begin{tabular}{|c|c|c|c|c|}
\hline
 & $[r1,r2,r3,r4,r5]$ & $Poincar\acute{e}\,\, series$ \\\hline
311 & $[1,1,3,4,4]$ & $\displaystyle{-\frac { \left( {t}^{2}+t+1 \right)
\left( {t}^{3
}+1 \right)  \left( {t}^{3}+{t}^{2}+t+1 \right) }{{t}^{7}+{t}^{6}+2\,{
t}^{5}+2\,{t}^{4}+2\,{t}^{3}+2\,t-1}}$\\\hline
312 & $[1,1,4,4,1]$ & $\displaystyle{-\frac { \left( {t}^{2}+t+1 \right)
\left( {t}^{3
}+{t}^{2}+t+1 \right)  \left( t+1 \right) }{{t}^{5}+2\,{t}^{4}+2\,{t}^
{3}+2\,{t}^{2}+t-1}}$\\\hline
313 & $[1,1,4,4,2]$ & $\displaystyle{-\frac { \left( {t}^{2}+t+1 \right)
\left( {t}^{3
}+{t}^{2}+t+1 \right)  \left( t+1 \right) }{{t}^{5}+2\,{t}^{4}+3\,{t}^
{3}+2\,{t}^{2}+t-1}}$\\\hline
314 & $[1,1,4,4,3]$ & $\displaystyle{-\frac { \left( {t}^{5}+{t}^{4}+{t}^
{3}+{t}^{2}+t+
1 \right)  \left( {t}^{3}+{t}^{2}+t+1 \right) }{{t}^{7}+{t}^{6}+2\,{t}
^{5}+2\,{t}^{4}+2\,{t}^{3}+2\,t-1}}$\\\hline
315 & $[1,1,4,4,4]$ & $\displaystyle{-\frac { \left( {t}^{2}+t+1 \right)
\left( {t}^{3
}+{t}^{2}+t+1 \right)  \left( t+1 \right) }{{t}^{6}+2\,{t}^{5}+3\,{t}^
{4}+3\,{t}^{3}+2\,{t}^{2}+t-1}}$\\\hline
316 & $[1,2,1,2,1]$ & $\displaystyle{\frac { \left( {t}^{3}+{t}^{2}+t+1
\right)
 \left( {t}^{2}+t+1 \right)  \left( {t}^{3}+1 \right)  \left( t+1
 \right) }{{t}^{9}+{t}^{8}+{t}^{6}-{t}^{5}-{t}^{4}-{t}^{3}-2\,{t}^{2}-
t+1}}$\\\hline
317 & $[1,2,1,2,2]$ & $\displaystyle{\frac { \left( {t}^{3}+{t}^{2}+t+1
\right)
 \left( {t}^{2}+t+1 \right)  \left( {t}^{3}+1 \right)  \left( t+1
 \right) }{{t}^{9}+{t}^{8}-{t}^{5}-{t}^{4}-2\,{t}^{3}-2\,{t}^{2}-t+1}}
]$\\
\hline

\end{tabular}

\begin{tabular}{|c|c|c|c|c|}
\hline
 & $[r1,r2,r3,r4,r5]$ & $Poincar\acute{e}\,\, series$ \\\hline

331 & $[1,2,2,1,4]$ & $\displaystyle{-\frac { \left( {t}^{3}+1 \right)
\left( {t}^{3}+
{t}^{2}+t+1 \right)  \left( {t}^{2}+t+1 \right)  \left( t+1 \right) }{
{t}^{7}+{t}^{6}+2\,{t}^{5}+2\,{t}^{4}+2\,{t}^{3}+2\,{t}^{2}+t-1}}$\\\hline
332 & $[1,2,2,2,1]$ & $\displaystyle{\frac { \left( {t}^{3}+1 \right)
\left( {t}^{3}+{
t}^{2}+t+1 \right)  \left( {t}^{2}+t+1 \right)  \left( t+1 \right) }{2
\,{t}^{9}+{t}^{8}-{t}^{5}-{t}^{4}-2\,{t}^{3}-2\,{t}^{2}-t+1}}$\\\hline
333 & $[1,2,2,2,2]$ & $\displaystyle{\frac { \left( {t}^{3}+1 \right)
\left( {t}^{3}+{
t}^{2}+t+1 \right)  \left( {t}^{2}+t+1 \right)  \left( t+1 \right) }{2
\,{t}^{9}+{t}^{8}-{t}^{6}-{t}^{5}-{t}^{4}-3\,{t}^{3}-2\,{t}^{2}-t+1}}$\\
\hline
334 & $[1,2,2,2,3]$ & $\displaystyle{\frac { \left( {t}^{3}+{t}^{2}+t+1
\right)
 \left( {t}^{5}+{t}^{4}+{t}^{3}+{t}^{2}+t+1 \right)  \left( t+1
 \right) }{2\,{t}^{9}+{t}^{8}-{t}^{6}-2\,{t}^{5}-2\,{t}^{4}-3\,{t}^{3}
-2\,{t}^{2}-t+1}}$\\\hline
335 & $[1,2,2,2,4]$ & $\displaystyle{\frac { \left( {t}^{2}+t+1 \right)
\left( {t}^{3}
+{t}^{2}+t+1 \right)  \left( {t}^{3}+1 \right)  \left( t+1 \right) }{{
t}^{9}-{t}^{7}-2\,{t}^{6}-2\,{t}^{5}-2\,{t}^{4}-3\,{t}^{3}-2\,{t}^{2}-
t+1}}$\\\hline
336 & $[1,2,2,3,1]$ & $\displaystyle{\frac { \left( {t}^{5}+{t}^{4}+{t}^
{3}+{t}^{2}+t+1
 \right)  \left( {t}^{3}+{t}^{2}+t+1 \right)  \left( t+1 \right) }{2\,
{t}^{9}+{t}^{8}-2\,{t}^{5}-2\,{t}^{4}-2\,{t}^{3}-2\,{t}^{2}-t+1}}$\\\hline
337 & $[1,2,2,3,2]$ & $\displaystyle{\frac { \left( {t}^{5}+{t}^{4}+{t}^
{3}+{t}^{2}+t+1
 \right)  \left( {t}^{3}+{t}^{2}+t+1 \right)  \left( t+1 \right) }{2\,
{t}^{9}+{t}^{8}-{t}^{6}-2\,{t}^{5}-2\,{t}^{4}-3\,{t}^{3}-2\,{t}^{2}-t+
1}}$\\
\hline

\end{tabular}

\begin{tabular}{|c|c|c|c|c|}
\hline
 & $[r1,r2,r3,r4,r5]$ & $Poincar\acute{e}\,\, series$ \\\hline

351 & $[1,2,3,2,4]$ & $\displaystyle{\frac { \left( {t}^{3}+1 \right)
\left( {t}^{2}+t
+1 \right)  \left( {t}^{3}+{t}^{2}+t+1 \right)  \left( t+1 \right) }{{
t}^{9}-{t}^{7}-2\,{t}^{6}-3\,{t}^{5}-3\,{t}^{4}-3\,{t}^{3}-2\,{t}^{2}-
t+1}}$\\\hline
352 & $[1,2,3,3,1]$ & $\displaystyle{\frac { \left( {t}^{3}+1 \right)
\left( {t}^{2}+t
+1 \right)  \left( {t}^{3}+{t}^{2}+t+1 \right)  \left( t+1 \right) }{2
\,{t}^{9}+{t}^{8}-3\,{t}^{5}-3\,{t}^{4}-2\,{t}^{3}-2\,{t}^{2}-t+1}}$\\\hline
353 & $[1,2,3,3,2]$ & $\displaystyle{\frac { \left( {t}^{3}+{t}^{2}+t+1
\right)
 \left( {t}^{5}+{t}^{4}+{t}^{3}+{t}^{2}+t+1 \right)  \left( t+1
 \right) }{2\,{t}^{9}+{t}^{8}-{t}^{6}-3\,{t}^{5}-3\,{t}^{4}-3\,{t}^{3}
-2\,{t}^{2}-t+1}}$\\\hline
354 & $[1,2,3,3,3]$ & $\displaystyle{\frac { \left( {t}^{3}+1 \right)
\left( {t}^{2}+t
+1 \right)  \left( {t}^{3}+{t}^{2}+t+1 \right)  \left( t+1 \right) }{2
\,{t}^{9}+{t}^{8}-{t}^{6}-4\,{t}^{5}-4\,{t}^{4}-3\,{t}^{3}-2\,{t}^{2}-
t+1}}$\\\hline
355 & $[1,2,3,3,4]$ & $\displaystyle{\frac { \left( {t}^{3}+1 \right)
\left( {t}^{2}+t
+1 \right)  \left( {t}^{3}+{t}^{2}+t+1 \right)  \left( t+1 \right) }{{
t}^{9}-{t}^{7}-2\,{t}^{6}-4\,{t}^{5}-4\,{t}^{4}-3\,{t}^{3}-2\,{t}^{2}-
t+1}}$\\\hline
356 & $[1,2,3,4,1]$ & $\displaystyle{\frac { \left( {t}^{3}+1 \right)
\left( {t}^{2}+t
+1 \right)  \left( {t}^{3}+{t}^{2}+t+1 \right)  \left( t+1 \right) }{{
t}^{9}-{t}^{7}-{t}^{6}-3\,{t}^{5}-3\,{t}^{4}-2\,{t}^{3}-2\,{t}^{2}-t+1
}}$\\\hline
357 & $[1,2,3,4,2]$ & $\displaystyle{\frac { \left( {t}^{3}+1 \right)
\left( {t}^{2}+t
+1 \right)  \left( {t}^{3}+{t}^{2}+t+1 \right)  \left( t+1 \right) }{{
t}^{9}-{t}^{7}-2\,{t}^{6}-3\,{t}^{5}-3\,{t}^{4}-3\,{t}^{3}-2\,{t}^{2}-
t+1}}$\\
\hline

\end{tabular}

\begin{tabular}{|c|c|c|c|c|}
\hline
 & $[r1,r2,r3,r4,r5]$ & $Poincar\acute{e}\,\, series$ \\\hline

371 & $[1,2,4,3,4]$ & $\displaystyle{-\frac { \left( {t}^{3}+{t}^{2}+t+1
\right)
 \left( {t}^{5}+{t}^{4}+{t}^{3}+{t}^{2}+t+1 \right)  \left( t+1
 \right) }{{t}^{8}+2\,{t}^{7}+3\,{t}^{6}+4\,{t}^{5}+4\,{t}^{4}+3\,{t}^
{3}+2\,{t}^{2}+t-1}}$\\\hline
372 & $[1,2,4,4,1]$ & $\displaystyle{-\frac { \left( {t}^{2}+t+1 \right)
\left( {t}^{3
}+{t}^{2}+t+1 \right)  \left( {t}^{3}+1 \right)  \left( t+1 \right) }{
{t}^{8}+2\,{t}^{7}+2\,{t}^{6}+3\,{t}^{5}+3\,{t}^{4}+2\,{t}^{3}+2\,{t}^
{2}+t-1}}$\\\hline
373 & $[1,2,4,4,2]$ & $\displaystyle{-\frac { \left( {t}^{2}+t+1 \right)
\left( {t}^{3
}+{t}^{2}+t+1 \right)  \left( {t}^{3}+1 \right)  \left( t+1 \right) }{
{t}^{8}+2\,{t}^{7}+3\,{t}^{6}+3\,{t}^{5}+3\,{t}^{4}+3\,{t}^{3}+2\,{t}^
{2}+t-1}}$\\\hline
374 & $[1,2,4,4,3]$ & $\displaystyle{-\frac { \left( {t}^{3}+{t}^{2}+t+1
\right)
 \left( {t}^{5}+{t}^{4}+{t}^{3}+{t}^{2}+t+1 \right)  \left( t+1
 \right) }{{t}^{8}+2\,{t}^{7}+3\,{t}^{6}+4\,{t}^{5}+4\,{t}^{4}+3\,{t}^
{3}+2\,{t}^{2}+t-1}}$\\\hline
375 & $[1,2,4,4,4]$ & $\displaystyle{-\frac { \left( {t}^{2}+t+1 \right)
\left( {t}^{3
}+{t}^{2}+t+1 \right)  \left( {t}^{3}+1 \right)  \left( t+1 \right) }{
{t}^{9}+2\,{t}^{8}+3\,{t}^{7}+4\,{t}^{6}+4\,{t}^{5}+4\,{t}^{4}+3\,{t}^
{3}+2\,{t}^{2}+t-1}}$\\\hline
376 & $[1,3,1,3,1]$ & $\displaystyle{\frac { \left( {t}^{5}+{t}^{4}+{t}^
{3}+{t}^{2}+t+1
 \right)  \left( t+1 \right) }{3\,{t}^{6}-2\,{t}^{5}-{t}^{4}+{t}^{3}-{
t}^{2}-2\,t+1}}$\\\hline
377 & $[1,3,1,3,2]$ & $\displaystyle{\frac { \left( {t}^{5}+{t}^{4}+{t}^
{3}+{t}^{2}+t+1
 \right)  \left( {t}^{3}+{t}^{2}+t+1 \right) }{3\,{t}^{8}-2\,{t}^{7}+2
\,{t}^{6}-2\,{t}^{5}-{t}^{4}-2\,{t}^{3}-2\,t+1}}$\\
\hline

\end{tabular}

\begin{tabular}{|c|c|c|c|c|}
\hline
 & $[r1,r2,r3,r4,r5]$ & $Poincar\acute{e}\,\, series$ \\\hline

391 & $[1,3,2,3,4]$ & $\displaystyle{\frac { \left( {t}^{5}+{t}^{4}+{t}^
{3}+{t}^{2}+t+1
 \right)  \left( {t}^{3}+{t}^{2}+t+1 \right) }{2\,{t}^{8}-2\,{t}^{7}+{
t}^{6}-3\,{t}^{5}-{t}^{4}-3\,{t}^{3}-2\,t+1}}$\\\hline
392 & $[1,3,2,4,1]$ & $\displaystyle{\frac { \left( {t}^{2}+t+1 \right)
\left( {t}^{3}
+{t}^{2}+t+1 \right)  \left( {t}^{3}+1 \right) }{2\,{t}^{8}-2\,{t}^{7}
+{t}^{6}-2\,{t}^{5}-{t}^{4}-2\,{t}^{3}-2\,t+1}}$\\\hline
393 & $[1,3,2,4,2]$ & $\displaystyle{\frac { \left( {t}^{2}+t+1 \right)
\left( {t}^{3}
+{t}^{2}+t+1 \right)  \left( {t}^{3}+1 \right) }{2\,{t}^{8}-2\,{t}^{7}
+{t}^{6}-3\,{t}^{5}-3\,{t}^{3}-2\,t+1}}$\\\hline
394 & $[1,3,2,4,3]$ & $\displaystyle{\frac { \left( {t}^{5}+{t}^{4}+{t}^
{3}+{t}^{2}+t+1
 \right)  \left( {t}^{3}+{t}^{2}+t+1 \right) }{2\,{t}^{8}-2\,{t}^{7}+{
t}^{6}-3\,{t}^{5}-{t}^{4}-3\,{t}^{3}-2\,t+1}}$\\\hline
395 & $[1,3,2,4,4]$ & $\displaystyle{\frac { \left( {t}^{2}+t+1 \right)
\left( {t}^{3}
+{t}^{2}+t+1 \right)  \left( {t}^{3}+1 \right) }{{t}^{8}-2\,{t}^{7}-3
\,{t}^{5}-{t}^{4}-3\,{t}^{3}-2\,t+1}}$\\\hline
396 & $[1,3,3,1,1]$ & $\displaystyle{\frac { \left( {t}^{2}+t+1 \right)
\left( {t}^{3}
+1 \right)  \left( t+1 \right) }{3\,{t}^{6}-2\,{t}^{5}-{t}^{4}+{t}^{3}
-{t}^{2}-2\,t+1}}$\\\hline
397 & $[1,3,3,1,2]$ & $\displaystyle{\frac { \left( {t}^{2}+t+1 \right)
\left( {t}^{3}
+{t}^{2}+t+1 \right)  \left( {t}^{3}+1 \right) }{3\,{t}^{8}-2\,{t}^{7}
+2\,{t}^{6}-2\,{t}^{5}-{t}^{4}-2\,{t}^{3}-2\,t+1}}$\\
\hline

\end{tabular}

\begin{tabular}{|c|c|c|c|c|}
\hline
 & $[r1,r2,r3,r4,r5]$ & $Poincar\acute{e}\,\, series$ \\\hline

411 & $[1,3,3,4,4]$ & $\displaystyle{\frac { \left( {t}^{2}+t+1 \right)

\left( {t}^{3}
+1 \right)  \left( t+1 \right) }{{t}^{6}-2\,{t}^{5}-{t}^{4}-{t}^{3}-{t
}^{2}-2\,t+1}}$\\\hline
412 & $[1,3,4,1,1]$ & $\displaystyle{\frac { \left( {t}^{2}+t+1 \right)

\left( {t}^{3}
+1 \right)  \left( t+1 \right) }{2\,{t}^{6}-2\,{t}^{5}-{t}^{4}+{t}^{3}
-{t}^{2}-2\,t+1}}$\\\hline
413 & $[1,3,4,1,2]$ & $\displaystyle{\frac { \left( {t}^{2}+t+1 \right)

\left( {t}^{3}
+{t}^{2}+t+1 \right)  \left( {t}^{3}+1 \right) }{2\,{t}^{8}-2\,{t}^{7}
+{t}^{6}-2\,{t}^{5}-{t}^{4}-2\,{t}^{3}-2\,t+1}}$\\\hline
414 & $[1,3,4,1,3]$ & $\displaystyle{\frac { \left( {t}^{2}+t+1 \right)

\left( {t}^{3}
+1 \right)  \left( t+1 \right) }{2\,{t}^{6}-2\,{t}^{5}-{t}^{4}-{t}^{2}
-2\,t+1}}$\\\hline
415 & $[1,3,4,1,4]$ & $\displaystyle{\frac { \left( {t}^{2}+t+1 \right)

\left( {t}^{3}
+1 \right)  \left( t+1 \right) }{{t}^{6}-2\,{t}^{5}-{t}^{4}-{t}^{2}-2
\,t+1}}$\\\hline
416 & $[1,3,4,2,1]$ & $\displaystyle{\frac { \left( {t}^{2}+t+1 \right)

\left( {t}^{3}
+{t}^{2}+t+1 \right)  \left( {t}^{3}+1 \right) }{2\,{t}^{8}-2\,{t}^{7}
+{t}^{6}-2\,{t}^{5}-{t}^{4}-2\,{t}^{3}-2\,t+1}}$\\\hline
417 & $[1,3,4,2,2]$ & $\displaystyle{\frac { \left( {t}^{2}+t+1 \right)

\left( {t}^{3}
+{t}^{2}+t+1 \right)  \left( {t}^{3}+1 \right) }{2\,{t}^{8}-2\,{t}^{7}
+{t}^{6}-3\,{t}^{5}-3\,{t}^{3}-2\,t+1}}$\\
\hline

\end{tabular}

\begin{tabular}{|c|c|c|c|c|}
\hline
 & $[r1,r2,r3,r4,r5]$ & $Poincar\acute{e}\,\, series$ \\\hline

431 & $[1,4,1,4,4]$ & $\displaystyle{-\frac { \left( {t}^{2}+t+1 \right)

\left( t+1
 \right) ^{2}}{2\,{t}^{3}+3\,{t}^{2}+t-1}}$\\\hline
432 & $[1,4,2,2,1]$ & $\displaystyle{\frac { \left( {t}^{3}+{t}^{2}+t+1

\right)
 \left( {t}^{2}+t+1 \right)  \left( t+1 \right) }{2\,{t}^{6}-{t}^{4}-3
\,{t}^{3}-2\,{t}^{2}-t+1}}$\\\hline
433 & $[1,4,2,2,2]$ & $\displaystyle{\frac { \left( {t}^{3}+{t}^{2}+t+1

\right)
 \left( {t}^{2}+t+1 \right)  \left( t+1 \right) }{2\,{t}^{6}-{t}^{4}-4
\,{t}^{3}-2\,{t}^{2}-t+1}}$\\\hline
434 & $[1,4,2,2,3]$ & $\displaystyle{\frac { \left( {t}^{5}+{t}^{4}+{t}^

{3}+{t}^{2}+t+1
 \right)  \left( {t}^{3}+{t}^{2}+t+1 \right) }{2\,{t}^{8}-2\,{t}^{7}+{
t}^{6}-3\,{t}^{5}-3\,{t}^{3}-2\,t+1}}$\\\hline
435 & $[1,4,2,2,4]$ & $\displaystyle{\frac { \left( {t}^{3}+{t}^{2}+t+1

\right)
 \left( {t}^{2}+t+1 \right)  \left( t+1 \right) }{{t}^{6}-{t}^{5}-2\,{
t}^{4}-4\,{t}^{3}-2\,{t}^{2}-t+1}}$\\\hline
436 & $[1,4,2,3,1]$ & $\displaystyle{\frac { \left( {t}^{5}+{t}^{4}+{t}^

{3}+{t}^{2}+t+1
 \right)  \left( {t}^{3}+{t}^{2}+t+1 \right) }{2\,{t}^{8}-2\,{t}^{7}+{
t}^{6}-2\,{t}^{5}-{t}^{4}-2\,{t}^{3}-2\,t+1}}$\\\hline
437 & $[1,4,2,3,2]$ & $\displaystyle{\frac { \left( {t}^{3}+{t}^{2}+t+1

\right)
 \left( {t}^{5}+{t}^{4}+{t}^{3}+{t}^{2}+t+1 \right) }{2\,{t}^{8}-2\,{t
}^{7}+{t}^{6}-3\,{t}^{5}-3\,{t}^{3}-2\,t+1}}$\\
\hline

\end{tabular}

\begin{tabular}{|c|c|c|c|c|}
\hline
 & $[r1,r2,r3,r4,r5]$ & $Poincar\acute{e}\,\, series$ \\\hline

451 & $[1,4,3,3,4]$ & $\displaystyle{\frac { \left( {t}^{5}+{t}^{4}+{t}^

{3}+{t}^{2}+t+1
 \right)  \left( t+1 \right) }{{t}^{6}-2\,{t}^{5}-{t}^{4}-{t}^{3}-{t}^
{2}-2\,t+1}}$\\\hline
452 & $[1,4,3,4,1]$ & $\displaystyle{\frac { \left( {t}^{5}+{t}^{4}+{t}^

{3}+{t}^{2}+t+1
 \right)  \left( t+1 \right) }{{t}^{6}-2\,{t}^{5}-{t}^{4}-{t}^{2}-2\,t
+1}}$\\\hline
453 & $[1,4,3,4,2]$ & $\displaystyle{\frac { \left( {t}^{5}+{t}^{4}+{t}^

{3}+{t}^{2}+t+1
 \right)  \left( {t}^{3}+{t}^{2}+t+1 \right) }{{t}^{8}-2\,{t}^{7}-3\,{
t}^{5}-{t}^{4}-3\,{t}^{3}-2\,t+1}}$\\\hline
454 & $[1,4,3,4,3]$ & $\displaystyle{\frac { \left( {t}^{5}+{t}^{4}+{t}^

{3}+{t}^{2}+t+1
 \right)  \left( t+1 \right) }{{t}^{6}-2\,{t}^{5}-{t}^{4}-{t}^{3}-{t}^
{2}-2\,t+1}}$\\\hline
455 & $[1,4,3,4,4]$ & $\displaystyle{-\frac { \left( {t}^{5}+{t}^{4}+{t}^

{3}+{t}^{2}+t+
1 \right)  \left( t+1 \right) }{2\,{t}^{5}+{t}^{4}+{t}^{3}+{t}^{2}+2\,
t-1}}$\\\hline
456 & $[1,4,4,1,1]$ & $\displaystyle{\frac { \left( {t}^{2}+t+1 \right)

\left( t+1
 \right) ^{2}}{{t}^{4}-{t}^{3}-3\,{t}^{2}-t+1}}$\\\hline
457 & $[1,4,4,1,2]$ & $\displaystyle{\frac { \left( {t}^{3}+{t}^{2}+t+1

\right)
 \left( {t}^{2}+t+1 \right)  \left( t+1 \right) }{{t}^{6}-{t}^{5}-2\,{
t}^{4}-3\,{t}^{3}-2\,{t}^{2}-t+1}}$\\
\hline

\end{tabular}

\begin{tabular}{|c|c|c|c|c|}
\hline
 & $[r1,r2,r3,r4,r5]$ & $Poincar\acute{e}\,\, series$ \\\hline

471 & $[1,4,4,4,4]$ & $\displaystyle{-\frac { \left( {t}^{2}+t+1 \right)

\left( t+1
 \right) ^{2}}{{t}^{4}+3\,{t}^{3}+3\,{t}^{2}+t-1}}$\\\hline
472 & $[2,2,2,2,1]$ & $\displaystyle{\frac { \left( {t}^{3}+{t}^{2}+t+1

\right)
 \left( {t}^{2}+t+1 \right)  \left( t+1 \right) }{3\,{t}^{6}+{t}^{5}-4
\,{t}^{3}-2\,{t}^{2}-t+1}}$\\\hline
473 & $[2,2,2,2,2]$ & $\displaystyle{\frac { \left( {t}^{3}+{t}^{2}+t+1

\right)
 \left( t+1 \right) }{3\,{t}^{4}-2\,{t}^{3}-{t}^{2}-2\,t+1}}$\\\hline
474 & $[2,2,2,2,3]$ & $\displaystyle{\frac { \left( {t}^{5}+{t}^{4}+{t}^

{3}+{t}^{2}+t+1
 \right)  \left( {t}^{3}+{t}^{2}+t+1 \right) }{3\,{t}^{8}-2\,{t}^{7}+2
\,{t}^{6}-4\,{t}^{5}+2\,{t}^{4}-4\,{t}^{3}-2\,t+1}}$\\\hline
475 & $[2,2,2,2,4]$ & $\displaystyle{\frac { \left( {t}^{3}+{t}^{2}+t+1

\right)
 \left( t+1 \right) }{2\,{t}^{4}-2\,{t}^{3}-{t}^{2}-2\,t+1}}$\\\hline
476 & $[2,2,2,3,1]$ & $\displaystyle{\frac { \left( {t}^{5}+{t}^{4}+{t}^

{3}+{t}^{2}+t+1
 \right)  \left( {t}^{3}+{t}^{2}+t+1 \right) }{3\,{t}^{8}-2\,{t}^{7}+2
\,{t}^{6}-3\,{t}^{5}+{t}^{4}-3\,{t}^{3}-2\,t+1}}$\\\hline
477 & $[2,2,2,3,2]$ & $\displaystyle{\frac { \left( {t}^{5}+{t}^{4}+{t}^

{3}+{t}^{2}+t+1
 \right)  \left( {t}^{3}+{t}^{2}+t+1 \right) }{3\,{t}^{8}-2\,{t}^{7}+2
\,{t}^{6}-4\,{t}^{5}+2\,{t}^{4}-4\,{t}^{3}-2\,t+1}}$\\
\hline

\end{tabular}

\begin{tabular}{|c|c|c|c|c|}
\hline
 & $[r1,r2,r3,r4,r5]$ & $Poincar\acute{e}\,\, series$ \\\hline

491 & $[2,2,3,4,4]$ & $\displaystyle{\frac { \left( {t}^{5}+{t}^{4}+{t}^

{3}+{t}^{2}+t+1
 \right)  \left( {t}^{3}+{t}^{2}+t+1 \right) }{{t}^{8}-2\,{t}^{7}-4\,{
t}^{5}-4\,{t}^{3}-2\,t+1}}$\\\hline
492 & $[2,2,4,4,1]$ & $\displaystyle{\frac { \left( {t}^{3}+{t}^{2}+t+1

\right)
 \left( {t}^{2}+t+1 \right)  \left( t+1 \right) }{{t}^{6}-{t}^{5}-2\,{
t}^{4}-4\,{t}^{3}-2\,{t}^{2}-t+1}}$\\\hline
493 & $[2,2,4,4,2]$ & $\displaystyle{\frac { \left( {t}^{3}+{t}^{2}+t+1

\right)
 \left( t+1 \right) }{{t}^{4}-2\,{t}^{3}-{t}^{2}-2\,t+1}}$\\\hline
494 & $[2,2,4,4,3]$ & $\displaystyle{\frac { \left( {t}^{5}+{t}^{4}+{t}^

{3}+{t}^{2}+t+1
 \right)  \left( {t}^{3}+{t}^{2}+t+1 \right) }{{t}^{8}-2\,{t}^{7}-4\,{
t}^{5}-4\,{t}^{3}-2\,t+1}}$\\\hline
495 & $[2,2,4,4,4]$ & $\displaystyle{-\frac { \left( {t}^{3}+{t}^{2}+t+1

\right)
 \left( t+1 \right) }{2\,{t}^{3}+{t}^{2}+2\,t-1}}$\\\hline
496 & $[2,3,2,3,1]$ & $\displaystyle{\frac { \left( {t}^{5}+{t}^{4}+{t}^

{3}+{t}^{2}+t+1
 \right)  \left( {t}^{3}+{t}^{2}+t+1 \right) }{3\,{t}^{8}-2\,{t}^{7}+2
\,{t}^{6}-3\,{t}^{5}-3\,{t}^{3}-2\,t+1}}$\\\hline
497 & $[2,3,2,3,2]$ & $\displaystyle{\frac { \left( {t}^{5}+{t}^{4}+{t}^

{3}+{t}^{2}+t+1
 \right)  \left( {t}^{3}+{t}^{2}+t+1 \right) }{3\,{t}^{8}-2\,{t}^{7}+2
\,{t}^{6}-4\,{t}^{5}+{t}^{4}-4\,{t}^{3}-2\,t+1}}$\\
\hline

\end{tabular}

\begin{tabular}{|c|c|c|c|c|}
\hline
 & $[r1,r2,r3,r4,r5]$ & $Poincar\acute{e}\,\, series$ \\\hline

511 & $[2,3,3,3,4]$ & $\displaystyle{\frac { \left( {t}^{5}+{t}^{4}+{t}^

{3}+{t}^{2}+t+1
 \right)  \left( {t}^{3}+{t}^{2}+t+1 \right) }{2\,{t}^{8}-2\,{t}^{7}+{
t}^{6}-4\,{t}^{5}-{t}^{4}-4\,{t}^{3}-2\,t+1}}$\\\hline
512 & $[2,3,3,4,1]$ & $\displaystyle{\frac { \left( {t}^{3}+1 \right)

\left( {t}^{2}+t
+1 \right)  \left( {t}^{3}+{t}^{2}+t+1 \right) }{2\,{t}^{8}-2\,{t}^{7}
+{t}^{6}-3\,{t}^{5}-{t}^{4}-3\,{t}^{3}-2\,t+1}}$\\\hline
513 & $[2,3,3,4,2]$ & $\displaystyle{\frac { \left( {t}^{5}+{t}^{4}+{t}^

{3}+{t}^{2}+t+1
 \right)  \left( {t}^{3}+{t}^{2}+t+1 \right) }{2\,{t}^{8}-2\,{t}^{7}+{
t}^{6}-4\,{t}^{5}-4\,{t}^{3}-2\,t+1}}$\\\hline
514 & $[2,3,3,4,3]$ & $\displaystyle{\frac { \left( {t}^{5}+{t}^{4}+{t}^

{3}+{t}^{2}+t+1
 \right)  \left( {t}^{3}+{t}^{2}+t+1 \right) }{2\,{t}^{8}-2\,{t}^{7}+{
t}^{6}-4\,{t}^{5}-{t}^{4}-4\,{t}^{3}-2\,t+1}}$\\\hline
515 & $[2,3,3,4,4]$ & $\displaystyle{\frac { \left( {t}^{5}+{t}^{4}+{t}^

{3}+{t}^{2}+t+1
 \right)  \left( {t}^{3}+{t}^{2}+t+1 \right) }{{t}^{8}-2\,{t}^{7}-4\,{
t}^{5}-{t}^{4}-4\,{t}^{3}-2\,t+1}}$\\\hline
516 & $[2,3,4,2,1]$ & $\displaystyle{\frac { \left( {t}^{2}+t+1 \right)

\left( {t}^{3}
+{t}^{2}+t+1 \right)  \left( {t}^{3}+1 \right) }{2\,{t}^{8}-2\,{t}^{7}
+{t}^{6}-3\,{t}^{5}-3\,{t}^{3}-2\,t+1}}$\\\hline
517 & $[2,3,4,2,2]$ & $\displaystyle{\frac { \left( {t}^{5}+{t}^{4}+{t}^

{3}+{t}^{2}+t+1
 \right)  \left( {t}^{3}+{t}^{2}+t+1 \right) }{2\,{t}^{8}-2\,{t}^{7}+{
t}^{6}-4\,{t}^{5}+{t}^{4}-4\,{t}^{3}-2\,t+1}}$\\
\hline

\end{tabular}

\begin{tabular}{|c|c|c|c|c|}
\hline
 & $[r1,r2,r3,r4,r5]$ & $Poincar\acute{e}\,\, series$ \\\hline

531 & $[2,4,2,4,4]$ & $\displaystyle{-\frac { \left( {t}^{3}+{t}^{2}+t+1

\right)
 \left( t+1 \right) }{2\,{t}^{3}+{t}^{2}+2\,t-1}}$\\\hline
532 & $[2,4,3,3,1]$ & $\displaystyle{\frac { \left( {t}^{5}+{t}^{4}+{t}^

{3}+{t}^{2}+t+1
 \right)  \left( {t}^{3}+{t}^{2}+t+1 \right) }{2\,{t}^{8}-2\,{t}^{7}+{
t}^{6}-3\,{t}^{5}-{t}^{4}-3\,{t}^{3}-2\,t+1}}$\\\hline
533 & $[2,4,3,3,2]$ & $\displaystyle{\frac { \left( {t}^{5}+{t}^{4}+{t}^

{3}+{t}^{2}+t+1
 \right)  \left( {t}^{3}+{t}^{2}+t+1 \right) }{2\,{t}^{8}-2\,{t}^{7}+{
t}^{6}-4\,{t}^{5}-4\,{t}^{3}-2\,t+1}}$\\\hline
534 & $[2,4,3,3,3]$ & $\displaystyle{\frac { \left( {t}^{5}+{t}^{4}+{t}^

{3}+{t}^{2}+t+1
 \right)  \left( {t}^{3}+{t}^{2}+t+1 \right) }{2\,{t}^{8}-2\,{t}^{7}+{
t}^{6}-4\,{t}^{5}-{t}^{4}-4\,{t}^{3}-2\,t+1}}$\\\hline
535 & $[2,4,3,3,4]$ & $\displaystyle{\frac { \left( {t}^{5}+{t}^{4}+{t}^

{3}+{t}^{2}+t+1
 \right)  \left( {t}^{3}+{t}^{2}+t+1 \right) }{{t}^{8}-2\,{t}^{7}-4\,{
t}^{5}-{t}^{4}-4\,{t}^{3}-2\,t+1}}$\\\hline
536 & $[2,4,3,4,1]$ & $\displaystyle{\frac { \left( {t}^{2}+t+1 \right)

\left( {t}^{3}
+1 \right)  \left( {t}^{3}+{t}^{2}+t+1 \right) }{{t}^{8}-2\,{t}^{7}-3
\,{t}^{5}-{t}^{4}-3\,{t}^{3}-2\,t+1}}$\\\hline
537 & $[2,4,3,4,2]$ & $\displaystyle{\frac { \left( {t}^{5}+{t}^{4}+{t}^

{3}+{t}^{2}+t+1
 \right)  \left( {t}^{3}+{t}^{2}+t+1 \right) }{{t}^{8}-2\,{t}^{7}-4\,{
t}^{5}-4\,{t}^{3}-2\,t+1}}$\\
\hline

\end{tabular}

\begin{tabular}{|c|c|c|c|c|}
\hline
 & $[r1,r2,r3,r4,r5]$ & $Poincar\acute{e}\,\, series$ \\\hline

551 & $[2,4,4,4,4]$ & $\displaystyle{-\frac { \left( {t}^{3}+{t}^{2}+t+1
\right)
 \left( t+1 \right) }{{t}^{4}+2\,{t}^{3}+{t}^{2}+2\,t-1}}$\\\hline
552 & $[3,3,3,3,1]$ & $\displaystyle{\frac { \left( {t}^{2}+t+1 \right)
\left( {t}^{3}
+1 \right)  \left( t+1 \right) }{3\,{t}^{6}-2\,{t}^{5}-{t}^{4}-{t}^{3}
-{t}^{2}-2\,t+1}}$\\\hline
553 & $[3,3,3,3,2]$ & $\displaystyle{\frac { \left( {t}^{5}+{t}^{4}+{t}^
{3}+{t}^{2}+t+1
 \right)  \left( {t}^{3}+{t}^{2}+t+1 \right) }{3\,{t}^{8}-2\,{t}^{7}+2
\,{t}^{6}-4\,{t}^{5}-{t}^{4}-4\,{t}^{3}-2\,t+1}}$\\\hline
554 & $[3,3,3,3,3]$ & $\displaystyle{\frac { \left( {t}^{5}+{t}^{4}+{t}^
{3}+{t}^{2}+t+1
 \right)  \left( t+1 \right) }{3\,{t}^{6}-2\,{t}^{5}-{t}^{4}-2\,{t}^{3
}-{t}^{2}-2\,t+1}}$\\\hline
555 & $[3,3,3,3,4]$ & $\displaystyle{\frac { \left( {t}^{5}+{t}^{4}+{t}^
{3}+{t}^{2}+t+1
 \right)  \left( t+1 \right) }{2\,{t}^{6}-2\,{t}^{5}-{t}^{4}-2\,{t}^{3
}-{t}^{2}-2\,t+1}}$\\\hline
556 & $[3,3,3,4,1]$ & $\displaystyle{\frac { \left( {t}^{2}+t+1 \right)
\left( {t}^{3}
+1 \right)  \left( t+1 \right) }{2\,{t}^{6}-2\,{t}^{5}-{t}^{4}-{t}^{3}
-{t}^{2}-2\,t+1}}$\\\hline
557 & $[3,3,3,4,2]$ & $\displaystyle{\frac { \left( {t}^{5}+{t}^{4}+{t}^
{3}+{t}^{2}+t+1
 \right)  \left( {t}^{3}+{t}^{2}+t+1 \right) }{2\,{t}^{8}-2\,{t}^{7}+{
t}^{6}-4\,{t}^{5}-{t}^{4}-4\,{t}^{3}-2\,t+1}}$\\
\hline

\end{tabular}

\begin{tabular}{|c|c|c|}
\hline
& $[r1,r2,r3,r4,r5,r6]$ & $Poincar\acute{e}\, series$ \\\hline
571 & $[3,4,4,3,4]$ & $\displaystyle{-\frac { \left( {t}^{5}+{t}^{4}+{t}^
{3}+{t}^{2}+t+
1 \right)  \left( t+1 \right) }{2\,{t}^{5}+{t}^{4}+2\,{t}^{3}+{t}^{2}+
2\,t-1}}$\\\hline
572 & $[3,4,4,4,1]$ & $\displaystyle{-\frac { \left( {t}^{2}+t+1 \right)
\left( {t}^{3
}+1 \right)  \left( t+1 \right) }{2\,{t}^{5}+{t}^{4}+{t}^{3}+{t}^{2}+2
\,t-1}}$\\\hline
573 & $[3,4,4,4,2]$ & $\displaystyle{-\frac { \left( {t}^{5}+{t}^{4}+{t}^
{3}+{t}^{2}+t+
1 \right)  \left( {t}^{3}+{t}^{2}+t+1 \right) }{2\,{t}^{7}+{t}^{6}+4\,
{t}^{5}+{t}^{4}+4\,{t}^{3}+2\,t-1}}$\\\hline
574 & $[3,4,4,4,3]$ & $\displaystyle{-\frac { \left( {t}^{5}+{t}^{4}+{t}^
{3}+{t}^{2}+t+
1 \right)  \left( t+1 \right) }{2\,{t}^{5}+{t}^{4}+2\,{t}^{3}+{t}^{2}+
2\,t-1}}$\\\hline
575 & $[3,4,4,4,4]$ & $\displaystyle{-\frac { \left( {t}^{5}+{t}^{4}+{t}^
{3}+{t}^{2}+t+
1 \right)  \left( t+1 \right) }{{t}^{6}+2\,{t}^{5}+{t}^{4}+2\,{t}^{3}+
{t}^{2}+2\,t-1}}$\\\hline
576 & $[4,4,4,4,1]$ & $\displaystyle{-\frac { \left( {t}^{2}+t+1 \right)
\left( t+1
 \right) ^{2}}{{t}^{4}+3\,{t}^{3}+3\,{t}^{2}+t-1}}$\\\hline
577 & $[4,4,4,4,2]$ & $\displaystyle{-\frac { \left( {t}^{3}+{t}^{2}+t+1
\right)
 \left( t+1 \right) }{{t}^{4}+2\,{t}^{3}+{t}^{2}+2\,t-1}}$\\\hline
578 & $[4,4,4,4,3]$ & $\displaystyle{-\frac { \left( {t}^{5}+{t}^{4}+{t}^
{3}+{t}^{2}+t+
1 \right)  \left( t+1 \right) }{{t}^{6}+2\,{t}^{5}+{t}^{4}+2\,{t}^{3}+
{t}^{2}+2\,t-1}}$\\\hline
579 & $[4,4,4,4,4]$ & $\displaystyle{-\frac { \left( t+1 \right) ^{2}}{2\,
{t}^{2}+2\,t-
1}}$\\\hline

\end{tabular}

\begin{center}
{\bf The Poincar\'{e} series of Cartan matrices with reduced Coxeter graph 4-VI}
\end{center}

\begin{center}

\end{center}

\begin{tabular}{|c|c|c|}
\hline
& $[r1,r2,r3,r4,r5,r6]$ & $Poincar\acute{e}\, series$ \\\hline

580 & $[1,1,1,1,1,1]$ & $\displaystyle{\frac { \left( {t}^{2}+t+1 \right)

\left( t+1
 \right) }{3\,{t}^{3}-2\,{t}^{2}-2\,t+1}}$\\\hline
581 & $[1,1,1,1,1,2]$ & $\displaystyle{\frac { \left( {t}^{3}+{t}^{2}+t+1

\right)
 \left( {t}^{2}+t+1 \right)  \left( t+1 \right) }{3\,{t}^{6}+{t}^{5}-{
t}^{4}-{t}^{3}-3\,{t}^{2}-t+1}}$\\\hline
582 & $[1,1,1,1,1,3]$ & $\displaystyle{\frac { \left( {t}^{5}+{t}^{4}+{t}^

{3}+{t}^{2}+t
+1 \right)  \left( t+1 \right) }{3\,{t}^{6}-2\,{t}^{5}-2\,{t}^{4}+3\,{
t}^{3}-2\,{t}^{2}-2\,t+1}}$\\\hline
583 & $[1,1,1,1,1,4]$ & $\displaystyle{\frac { \left( {t}^{2}+t+1 \right)

\left( t+1
 \right) }{2\,{t}^{3}-2\,{t}^{2}-2\,t+1}}$\\\hline
584 & $[1,1,1,1,2,2]$ & $\displaystyle{\frac { \left( {t}^{3}+{t}^{2}+t+1

\right)
 \left( {t}^{2}+t+1 \right)  \left( t+1 \right) }{3\,{t}^{6}+{t}^{5}-{
t}^{4}-2\,{t}^{3}-3\,{t}^{2}-t+1}}$\\\hline
585 & $[1,1,1,1,2,3]$ & $\displaystyle{\frac { \left( {t}^{5}+{t}^{4}+{t}^

{3}+{t}^{2}+t
+1 \right)  \left( {t}^{3}+{t}^{2}+t+1 \right) }{3\,{t}^{8}-2\,{t}^{7}
+{t}^{6}-3\,{t}^{4}-{t}^{2}-2\,t+1}}$\\\hline
586 & $[1,1,1,1,2,4]$ & $\displaystyle{\frac { \left( {t}^{3}+{t}^{2}+t+1

\right)
 \left( {t}^{2}+t+1 \right)  \left( t+1 \right) }{2\,{t}^{6}-2\,{t}^{4
}-2\,{t}^{3}-3\,{t}^{2}-t+1}}$\\
\hline

\end{tabular}%

\end{center}
\end{document}